\documentclass[12pt,letterpaper,reqno]{amsart}
\usepackage{amsfonts}
\usepackage{color}
\usepackage{epsfig}
\usepackage{fancyhdr}
\usepackage[latin1]{inputenc}
\usepackage{amsmath}
\usepackage{mathrsfs}
\usepackage{natbib}
\usepackage{graphicx}
\usepackage[colorlinks=true,linkcolor=red,citecolor=blue]{hyperref}
\usepackage{amssymb}
\usepackage{dsfont}
\usepackage{fancyhdr}
\usepackage{graphicx}

\newtheorem{theorem}{Theorem}

\newtheorem{corollary}[theorem]{Corollary}

\newtheorem{definition}[theorem]{Definition}
\newtheorem{example}[theorem]{Example}

\newtheorem{proposition}[theorem]{Proposition}
\newtheorem{remark}[theorem]{Remark}

\renewcommand{\cite}{\citet}
 
\newcommand{\E }{\mathbb{E}}

\newcommand{\R }{\mathbb{R}}

\begin{document}
\author{Michel Broniatowski$^{1}$, Jana Jure\v{c}ková$^2$ and Amor KEZIOU$%
^{3}$}
\address{$^{1}$ LPSM, Sorbonne-Université, Paris, France \\
michel.broniatowski@upmc.fr}
\address{$^{2}$The Czech Academy of Sciences, Institute of Information
Theory and Automation, Prague, Czech Republic\\
jureckova@karlin.mff.cuni.cz}
\address{$^{3}$Laboratoire de Mathématiques de Reims, France\\
amor.keziou@univ-reims.fr}
\date{April, 2019}
\title{Uniform minimum risk equivariant estimates for moment condition models%
}
\maketitle

\begin{abstract}
We consider semiparametric moment condition models invariant to
transformation groups. The parameter of interest is estimated by minimum
empirical divergence approach, introduced by \cite{Bronia_Keziou2012}. It is
shown that the minimum empirical divergence estimates, including the
empirical likelihood one, are equivariants. The minimum risk equivariant
estimate is then identified to be any one of the minimum empirical
divergence estimates minus its expectation conditionally to maximal
invariant statistic of the considered group of transformations. An
asymptotic approximation to the conditional expectation, is obtained, using
the result of \cite{Jurekova_Picek_2009_StatDec}. \vskip0.5cm \noindent 
\end{abstract}

\section{Introduction}

\noindent The semiparametric moment condition models are defined through
estimating equations 
\begin{equation*}
\mathbb{E}\left(f_j(\mathbf{X},\boldsymbol{\theta}_T)\right) = 0 \text{ for
all } j=1,\ldots,\ell,
\end{equation*}
where $\mathbb{E}(\cdot)$ denotes the mathematical expectation, $\boldsymbol{%
X}\in\mathbb{R}^m$ is a random vector, $\boldsymbol{\theta}%
_T\in\Theta\subseteq\mathbb{R}^d$ is the unknown true value of the parameter
of interest which is assumed to be unique, and $\mathbf{f}(\boldsymbol{x},%
\boldsymbol{\theta}):= (f_1(\boldsymbol{x},\boldsymbol{\theta}%
),\ldots,f_\ell(\boldsymbol{x},\boldsymbol{\theta}))^\top$ is some specified
measurable $\mathbb{R}^\ell$-valued function defined on $\mathbb{R}^m\times
\Theta$. Such models are popular in statistics and econometrics, see e.g.,
\noindent\cite{Qin-Lawless1994}, \cite{Haberman1984}, \cite{Sheehy1987}, 
\cite{McCullagh_Nelder1983}, \cite{Owen2001} and the references therein.
Denoting $P_{\boldsymbol{X}}(\cdot)$ the probability distribution of the
random vector $\boldsymbol{X}$, then the above estimating equations can be
written as 
\begin{equation*}
\int_{\mathbb{R}^m} \mathbf{f}(\boldsymbol{x},\boldsymbol{\theta}_T)\,dP_{%
\boldsymbol{X}}(\boldsymbol{x}) = \boldsymbol{0}.
\end{equation*}
Let $M^{1}$ be the collection of all signed finite measures (s.f.m.) $Q$ on
the Borel $\sigma$-field $(\mathbb{R}^m,\mathcal{B}(\mathbb{R}^m))$ such
that $Q(\mathbb{R}^m)=1$. The submodel $\mathcal{M}_{\boldsymbol{\theta}}$,
associated to a given value $\boldsymbol{\theta}\in\Theta$, consists of all
s.f.m.'s $Q\in M^1$ satisfying $\ell$ linear constraints induced by the
vector valued function $\mathbf{f}(\cdot,\boldsymbol{\theta}):=(f_1(\cdot,%
\boldsymbol{\theta}),\ldots,f_\ell(\cdot,\boldsymbol{\theta}))^\top$,
namely, 
\begin{eqnarray}
\mathcal{M}_{\boldsymbol{\theta}} & := & \left\{Q\in M^{1}\text{ such that }%
\int_{\mathbb{R}^m} f_j(\boldsymbol{x},\boldsymbol{\theta})\,dQ(\boldsymbol{x%
})=0,\,\forall j=1,\ldots,\ell\right\}  \notag \\
& = & \left\{Q\in M^{1}\text{ such that } \int_{\mathbb{R}^m} \mathbf{f}(%
\boldsymbol{x},\boldsymbol{\theta})\,dQ(\boldsymbol{x}) = \mathbf{0}\right\},
\notag
\end{eqnarray}
with $\ell > d$. The statistical model which we consider can be written as 
\begin{equation}  \label{the model}
\mathcal{M} := \bigcup_{\boldsymbol{\theta}\in\Theta}\mathcal{M}_{%
\boldsymbol{\theta}} := \bigcup_{\boldsymbol{\theta}\in\Theta} \left\{Q\in
M^{1}\text{ such that } \int_{\mathbb{R}^m} \mathbf{f}(\boldsymbol{x},%
\boldsymbol{\theta})\,dQ(\boldsymbol{x}) = \boldsymbol{0}\right\}.
\end{equation}
\noindent Let ${\boldsymbol{X}}_1,\ldots,{\boldsymbol{X}}_n$ be an i.i.d.
sample of the random vector $\boldsymbol{X}\in\mathbb{R}^m$ with unknown
probability distribution $P_{\boldsymbol{X}}(\cdot)$. The problems of
testing the model $\mathcal{H}_0 : P_{\boldsymbol{X}}\in \mathcal{M}$,
confidence region and point estimations of $\boldsymbol{\theta}_T$, have
been widely investigated in the literature. \cite{Hansen1982} considered
generalized method of moments (GMM) in order to estimate $\boldsymbol{\theta}%
_T$. \cite{Hansen_Healton_Yaron1996} introduced the continuous updating (CU)
estimate. Asymptotic confidence regions for the parameter $\boldsymbol{\theta%
}_T$ have been obtained by \cite{Owen1988} and \cite{Owen1990}, introducing
the empirical likelihood (EL) approach. It has been used, in the context of
model (\ref{the model}), by \cite{Qin-Lawless1994} and \cite{Imbens1997}
introducing the EL estimate for the parameter $\boldsymbol{\theta}_T$. The
recent literature in econometrics focusses on such models; \cite{Smith1997}, 
\cite{NeweySmith2004} provided a class of estimates called generalized
empirical likelihood (GEL) estimates which contains the EL and the CU ones.
Among other results pertaining to EL, \cite{NeweySmith2004} stated that EL
estimate enjoys asymptotic optimality properties in term of efficiency when
bias corrected among all GEL estimates including the GMM one. \cite%
{Bronia_Keziou2012} proposed a general approach through empirical
divergences and duality technique which includes the above methods in the
general context of signed finite measures under moment condition models (\ref%
{the model}). These approach allows the asymptotic study of the estimates
and associated test statistics both under the model and under
misspecification, leading to new results, in particular, for the EL
approach. Note that all the proposed estimates including the EL one are
generally biased, and that the problem of their finite sample efficiency, at
our knowledge, have not yet been studied.\newline

The aim of the present paper is to investigate the finite-sample optimality
property estimation in the context of semiparametric model (\ref{the model}%
). We will discuss the problem of constructing minimum risk equivariant
estimates (MRE) for the parameter $\boldsymbol{\theta}_T$, as well as the
problem of the numerical calculation of these estimates.\newline

We recall in the following lines, for the above estimation problem, the
notions of group transformations on the random vector space, model
invariance and the induced group of transformations on the parameter space,
loss invariance and equivariance estimation; we refer to the unpublished
preprint of \cite{Hoff_2012_Preprint} for an excellent presentation of the
above notions, and the book of \cite{Lehmann_Casella1998}.\newline

Let $\mathcal{G}$ be a collection, of one-to-one transformations from the
vector space $\mathbb{R}^{m}$ in $\mathbb{R}^{m}$, which we assume to be a
``group'', in the sense that, it should be closed under both composition and
inversion, namely, 
\begin{equation*}
\forall g_1, g_2\in \mathcal{G},\, g_1 \circ g_2\in\mathcal{G} \quad \text{
and }\quad \forall g\in\mathcal{G},\, g^{-1}\in\mathcal{G}.
\end{equation*}
The group $\mathcal{G}$ can be extended to a group of transformations on the
sample space, $\R^{mn}$ onto $\R^{mn}$, which will be denoted $\mathcal{G}_n$%
, as follows 
\begin{equation*}  \label{group on the sample sapce}
\mathcal{G}_n:=\left\{\mathbf{g} : (\boldsymbol{x}_1,\ldots,\boldsymbol{x}%
_n)\in\mathbb{R}^{mn}\mapsto \mathbf{g}(\boldsymbol{x}_1,\ldots,\boldsymbol{x%
}_n):=(g(\boldsymbol{x}_1), \ldots,g(\boldsymbol{x}_n))\in\mathbb{R}^{mn};\,
g\in\mathcal{G}\right\}.
\end{equation*}

We will consider two kinds of transformation groups,

\begin{enumerate}
\item[$\bullet$ ] ``additive'' 
\begin{equation}  \label{additive group}
\mathcal{G} :=\left\{\boldsymbol{x}\in\R^m \mapsto \boldsymbol{x}+%
\boldsymbol{a} \in\R^m;\, \boldsymbol{a}\in S\right\},
\end{equation}
where $S$ is some subset of $\R^m$,

\item[or] 

\item[$\bullet$] ``multiplicative'' 
\begin{equation}  \label{multiplicative group}
\mathcal{G} :=\left\{\boldsymbol{x}\in\R^m \mapsto \text{diag}%
(\lambda_1,\ldots,\lambda_m) \boldsymbol{x} \in\R^m; \mathbf{\lambda}\in\R^m 
\text{ or } \mathbf{\lambda}\in {\R^*}^m\right\},
\end{equation}
where $\text{diag}(\lambda_1,\ldots,\lambda_m)$ is diagonal matrix, with
entries $\lambda_1,\ldots\lambda_m\in\R^*$ or $\lambda_1,\ldots\lambda_m\in\R%
^*_+$ with possibly some entries $\lambda_i$ equal to one.\newline
\end{enumerate}

We assume that the model $\mathcal{M}$ given in (\ref{the model}) is
invariant under the considered group of transformations $\mathcal{G}$, in
the sense that, 
\begin{equation}  \label{invariance condition}
\text{ for any random vector } \boldsymbol{X}, \text{ if } P_{\boldsymbol{X}%
}\in\mathcal{M}, \text{ then } P_{g(\boldsymbol{X})}\in\mathcal{M}, \forall
g\in\mathcal{G}. \\
\end{equation}

The induced group of transformations on the parameter space, $\Theta$ onto $%
\Theta$, denoted $\overline{\mathcal{G}}$ hereafter, will be defined as
follows. Let $g$ be any transformation belonging to $\mathcal{G}$, and
consider any random vector ${\boldsymbol{X}}$ such that $P_{\boldsymbol{X}%
}\in\mathcal{M}$. Then, by identifiability assumption, there exists a unique 
$\boldsymbol{\theta}\in\Theta$ such that $P_{\boldsymbol{X}}\in\mathcal{M}_{%
\boldsymbol{\theta}}$. By invariance assumption (\ref{invariance condition}%
), of the model $\mathcal{M}$ to the group $\mathcal{G}$, the distribution $%
P_{g(\mathbf{X})}$ belongs to $\mathcal{M}$. Therefore, there exists a
unique (by indentifiability) $\overline{\boldsymbol{\theta}}\in\Theta$ such
that $P_{g(\mathbf{X})}\in \mathcal{M}_{\overline{\boldsymbol{\theta}}}$.
Denote then by $\overline{g}$ the bijection induced by $g$ on the parameter
space $\Theta$ onto $\Theta$, defined by 
\begin{equation*}
\overline{g} : \boldsymbol{\theta}\in\Theta \mapsto \overline{g}(\boldsymbol{%
\theta}) := \overline{\boldsymbol{\theta}}\in\Theta.
\end{equation*}
The induced group on the parameter space, $\Theta$ onto $\Theta$, is then
defined to be 
\begin{equation*}  \label{the group on the parameter space}
\overline{\mathcal{G}} :=\left\{\overline{g} \text{ such that } g\in\mathcal{%
G}\right\}.
\end{equation*}

Two points $\boldsymbol{\theta}_1, \boldsymbol{\theta}_2\in\Theta$ are said
equivalent iff $\boldsymbol{\theta}_2 = \overline{g}(\boldsymbol{\theta}_1)$
for some $\overline{g}\in\overline{\mathcal{G}}.$ The orbit $\Theta(%
\boldsymbol{\theta}_0)$, of a point $\boldsymbol{\theta}_0\in\Theta$, is
defined to be the set of equivalent points: 
\begin{equation*}  \label{the orbit of a point}
\Theta(\boldsymbol{\theta}_0) := \left\{\overline{g}(\boldsymbol{\theta}%
_0);\, \overline{g}\in\overline{G}\right\}.
\end{equation*}
We will assume that there is only one orbit of $\Theta$, i.e., 
\begin{equation}  \label{transitive assumption}
\Theta(\boldsymbol{\theta}_0) =\Theta
\end{equation}
which means that the group of transformation $\overline{\mathcal{G}}$ is
rich enough allowing to go from any point in $\Theta$ to another via some
transformation $\overline{g}\in\overline{\mathcal{G}}$. In such case, the
group $\overline{\mathcal{G}}$ is said to be ``transitive'' over $\Theta$. 
\newline

We give here some examples for illustration. In all the examples below, we
can see that the group $\overline{\mathcal{G}}$ is transitive over $\Theta$.

\begin{example}
\label{example 1} Sometimes we have information relating the first and
second moments of a random variable $X$ (see e.g. \cite%
{GodambeAndThompson1989} and \cite{McCullagh_Nelder1983}). Let $%
X_1,\ldots,X_n$ be an i.i.d. sample of a random variable $X\in\mathbb{R}$
with mean $\mathbb{E}(X)=\theta_T$, and assume that $\mathbb{E}(X^2) =
h(\theta_T)$, where $h(\cdot)$ is a known function. Our aim is to estimate $%
\theta_T$. The information about the distribution $P_X$ of $X$ can be
expressed in the form of (\ref{the model}) by taking $\mathbf{f}%
(x,\theta):=\left(x-\theta,x^2-h(\theta)\right)^\top.$ If we take the
parameter space to be $\Theta=\mathbb{R}$, then it is straightforward to see
that the model $\mathcal{M}$ is invariant to the additive group of
transformations 
\begin{equation*}
\mathcal{G} :=\left\{g : x\in\mathbb{R}\mapsto g(x):= x+a; \, a\in \mathbb{R}%
\right\},
\end{equation*}
if $h(\theta) := \theta^2+c$ for some $c \geq 0$, and invariant to the
multiplicative group 
\begin{equation*}
\mathcal{G} := \left\{ g : x\in\mathbb{R}\mapsto g(x):= \lambda x; \,
\lambda\in \mathbb{R}^*_+\right\},
\end{equation*}
if $h(\theta) := c\theta^2$ for some $c > 0$.The induced groups on the
parameter space $\Theta :=\mathbb{R}$ are, respectively, 
\begin{equation*}
\overline{\mathcal{G}}=\left\{\overline{g} : \theta\in\mathbb{R}\mapsto 
\overline{g}(\theta) = \theta+a;\, a\in\mathbb{R}\right\}
\end{equation*}
and 
\begin{equation*}
\overline{\mathcal{G}}=\left\{\overline{g} : \theta\in\mathbb{R}\mapsto 
\overline{g}(\theta) = \lambda\theta;\, \lambda\in\mathbb{R}^*_+\right\}.
\end{equation*}
\end{example}

\begin{example}
\label{example 2} Let $\boldsymbol{X}_1:=(X_{1,1},X_{2,1})^\top,\ldots,%
\boldsymbol{X}_n:=(X_{1,n},X_{2,n})^\top$ be an i.i.d. sample of a bivariate
random vector $\boldsymbol{X} := (X_1,X_2)^\top$ with $\mathbb{E}(X_1)=%
\mathbb{E}(X_2) = \theta_T$. In this case, we can take $\mathbf{f}(\mathbf{x}%
,\theta)=(x_1-\theta, x_2-\theta)^\top.$ If we consider $\Theta=\mathbb{R}$,
then the model $\mathcal{M}$ is invariant with respect to the groups 
\begin{equation*}
\mathcal{G} :=\left\{g : \boldsymbol{x}\in\mathbb{R}^2\mapsto g(\boldsymbol{x%
}):= \boldsymbol{x}+a\mathbf{1};\, a\in \mathbb{R}\right\}
\end{equation*}
or 
\begin{equation*}
\mathcal{G} :=\left\{ g : \boldsymbol{x}\in\mathbb{R}^2\mapsto g(\boldsymbol{%
x}) := \lambda \boldsymbol{x};\, \lambda\in \mathbb{R}^*_+\right\}.
\end{equation*}
The induced groups on $\Theta:=\mathbb{R}$ are, respectively, 
\begin{equation*}
\overline{\mathcal{G}}=\left\{\overline{g}: \theta\in\mathbb{R}\mapsto 
\overline{g}(\theta) = \theta+a;\, a\in\mathbb{R}\right\}
\end{equation*}
and 
\begin{equation*}
\overline{\mathcal{G}}=\left\{\overline{g}:\theta\in\mathbb{R}\mapsto 
\overline{g}(\theta) = \lambda\theta;\, \lambda\in\mathbb{R}^*_+\right\}.
\end{equation*}
A some what similar problem is when $\mathbb{E}(X_1)=c$ is known, and $%
\mathbb{E}(X_2)=\theta_T$ is to be estimated, by taking $\mathbf{f}(\mathbf{x%
},\theta)=(x_1-c,x_2-\theta)^\top.$ Such problems are common in survey
sampling (see e.g. \cite{KukAndMak1989} and \cite{ChenAndQin1993}). Taking $%
\Theta=\mathbb{R}$, the model $\mathcal{M}$ is then invariant with respect
to the groups 
\begin{equation*}
\mathcal{G} :=\left\{g : \boldsymbol{x}\in\mathbb{R}^2\mapsto g(\boldsymbol{x%
}):= \boldsymbol{x}+(0,a)^\top\in\R^2; \, a\in \mathbb{R}\right\}
\end{equation*}
or 
\begin{equation*}
\mathcal{G} :=\left\{ g : \boldsymbol{x}\in\mathbb{R}^2\mapsto g(\boldsymbol{%
x}):= (x_1,\lambda x_2)^\top\in\R^2;\, \lambda\in \mathbb{R}^*_+\right\}.
\end{equation*}
The induced groups on $\Theta:=\mathbb{R}$ are, respectively, 
\begin{equation*}
\overline{\mathcal{G}}=\left\{\overline{g}:\theta\in\mathbb{R}\mapsto 
\overline{g}(\theta) = \theta+a;\, a\in\mathbb{R}\right\}
\end{equation*}
and 
\begin{equation*}
\overline{\mathcal{G}}=\left\{\overline{g}:\theta\in\mathbb{R}\mapsto 
\overline{g}(\theta) = \lambda\theta;\, \lambda\in\mathbb{R}^*_+\right\}.
\end{equation*}
\end{example}

\begin{example}
\label{example 3} Let $X_1,\ldots,X_n$ be an i.i.d. sample of a random
variable $X$ with distribution $P_X$ such that $\mathbb{E}\left(f_i(X -
\theta_T)\right) = 0$, where $f_i:x\in\mathbb{R}\mapsto f_i(x) := \mathds{1}%
_{(a_i,b_i)}(x) -c_i$, $\forall i=1,\ldots,\ell$. The known intervals $%
(a_i,b_i)$ may be bounded or unbounded, and $c_1,\ldots,c_\ell$ are known
nonnegative numbers. The information about $P_X$ can be written under the
form of model (\ref{the model}) taking $\mathbf{f}(x,\theta):=(f_1(x-%
\theta),\ldots,f_\ell(x-\theta))^\top$ and $\theta\in\Theta :=\mathbb{R}$.
The model $\mathcal{M}$ in this case is invariant to the groups 
\begin{equation*}
\mathcal{G} :=\left\{g : x\in\mathbb{R}\mapsto g(x):= x+a;\, a\in \mathbb{R}%
\right\}
\end{equation*}
or 
\begin{equation*}
\mathcal{G} :=\left\{g : x\in\mathbb{R}\mapsto g(x):= \lambda x;\,
\lambda\in \mathbb{R}_+^*\right\},
\end{equation*}
and the induced groups on the parameter space $\Theta$ are, respectively, 
\begin{equation*}
\overline{\mathcal{G}} = \left\{\overline{g}:\theta\in\mathbb{R}\mapsto 
\overline{g}(\theta) = \theta+a;\, a\in\mathbb{R}\right\}
\end{equation*}
and 
\begin{equation*}
\overline{\mathcal{G}} = \left\{\overline{g}:\theta\in\mathbb{R}\mapsto 
\overline{g}(\theta) = \lambda\theta;\, \lambda\in\mathbb{R}_+^*\right\}.
\end{equation*}
\end{example}

\begin{example}
\label{example 4_1} Let $X_1,\ldots,X_n$ be an i.i.d. sample of a random
variable $X\in\mathbb{R}$ with continuous distribution $P_X$ such that $%
\mathbb{E}(X)=\theta_{1T}$, $\mathbb{E}(X^2)=1+\theta_{1T}^2$ and $\mathbb{E}%
\left(\mathds{1}_{]-\infty,\theta_{2T}]}(X)\right) = \alpha$, where $\alpha
\in ]0,1[$ is known and $\boldsymbol{\theta}_T:=(\theta_{1T},\theta_{2T})^%
\top$ is to be estimated. Note that $\theta_{2T}$ is the quantile of order $%
\alpha$ of the variable $X$, and that the variance of $X$ is assumed to be
known and equal to one. This problem can be written under the form of model (%
\ref{the model}) taking $\mathbf{f}(x,\boldsymbol{\theta}) :=
\left(x-\theta_1, x^2-1-\theta_1^2,\mathds{1}_{]-\infty,0]}(x-\theta)%
\right)^\top$ and $\boldsymbol{\theta}\in\Theta :=\mathbb{R}^2$. The model $%
\mathcal{M}$ in this case is invariant with respect to the additive group 
\begin{equation*}
\mathcal{G} :=\left\{g : x\in\mathbb{R}\mapsto g(x):= x+a;\, a\in \mathbb{R}%
\right\}
\end{equation*}
and the induced group on $\Theta$ is 
\begin{equation*}
\overline{\mathcal{G}}=\left\{\overline{g}:\boldsymbol{\theta}\in\mathbb{R}%
^2\mapsto \overline{g}(\boldsymbol{\theta}) = \boldsymbol{\theta}+a\mathbf{1}%
;\, a\in\mathbb{R}\right\}.
\end{equation*}
\end{example}

\begin{example}
\label{example 4_2} Let $X_1,\ldots,X_n$ be an i.i.d. sample of a random
variable $X\in\mathbb{R}$ with continuous distribution $P_X$ such that $%
\mathbb{E}(X)=0$ and $\mathbb{E}\left(\mathds{1}_{]-\infty,\theta_T]}(X)%
\right) = \alpha$, where $\alpha \in ]0,1[$ is known and $\theta_T$ is to be
estimated. Note that $\theta_T$ is the quantile of order $\alpha$ of the
variable $X$. This problem can be written under the form of model (\ref{the
model}) taking $\mathbf{f}(x,\theta) := \left(x,\mathds{1}%
_{]-\infty,0]}(x-\theta)\right)$ and $\theta\in\Theta :=\mathbb{R}$. The
model $\mathcal{M}$ in this case is invariant with respect to the
multiplicative group 
\begin{equation*}
\mathcal{G} :=\left\{g : x\in\mathbb{R}\mapsto g(x):= \lambda x;\, a\in 
\mathbb{R}_+^*\right\}
\end{equation*}
and the induced group on $\Theta$ is 
\begin{equation*}
\overline{\mathcal{G}}=\left\{\overline{g}:\theta\in\mathbb{R}\mapsto 
\overline{g}(\theta) = \lambda \theta;\, \lambda\in\mathbb{R}_+^*\right\}.
\end{equation*}
\end{example}

\begin{example}
\label{example 5} Let $\boldsymbol{X}_1,\ldots,\boldsymbol{X}_n$ be an
i.i.d. sample of a random vector $\boldsymbol{X}\in\mathbb{R}^m$ with
continuous distribution $P_{\boldsymbol{X}}$ such that $\mathbb{E}\left(%
\mathbf{h}(\boldsymbol{X}-\boldsymbol{\theta}_T)\right) = 0$, where $\mathbf{%
h} : \R^m\mapsto \R^\ell$ is some specified measurable function, and $%
\boldsymbol{\theta}_T\in\R^m$ is to be estimated. We can consider also the
case where some components of $\boldsymbol{\theta}_T$ are known and that the
other components are to be estimated. It is clear that the corresponding
model $\mathcal{M}$ defined in (\ref{the model}), taking $\mathbf{f}(%
\boldsymbol{x},\boldsymbol{\theta}):=\mathbf{h}(\boldsymbol{x}-\boldsymbol{%
\theta})$ and $\Theta = \R^m$, is invariant to the additive group 
\begin{equation*}
\mathcal{G} :=\left\{g : \boldsymbol{x}\in\R^m\mapsto g(\boldsymbol{x}) := 
\boldsymbol{x}+\boldsymbol{a}; \, \boldsymbol{a}\in \R^m\right\}
\end{equation*}
and the induced group on the parameter is 
\begin{equation*}
\overline{\mathcal{G}}=\left\{\overline{g}:\boldsymbol{\theta}\in\R^m\mapsto 
\overline{g}(\boldsymbol{\theta}) = \boldsymbol{\theta} + \boldsymbol{a};\, 
\boldsymbol{a}\in\R^m\right\}.
\end{equation*}
Likewise, if the data are such that $\mathbb{E}\left(\mathbf{h}\left(\frac{%
\boldsymbol{X}}{{\boldsymbol{\theta}}_T}\right)\right) = 0$, where $\mathbf{h%
} : \R^m\mapsto \R^\ell$ is some specified measurable function, and $%
\boldsymbol{\theta}_T\in\Theta\subset\R^m$ is to be estimated, then the
corresponding model $\mathcal{M}$, taking $\mathbf{f}(\boldsymbol{x},%
\boldsymbol{\theta}):=\mathbf{h}\left(\frac{\boldsymbol{x}}{\boldsymbol{%
\theta}}\right)$, is invariant to the multiplicative group 
\begin{equation*}
\mathcal{G} :=\left\{g : \boldsymbol{x}\in\R^m\mapsto g(\boldsymbol{x}) := 
\text{diag}(\lambda_1,\ldots,\lambda_m)\boldsymbol{x}; \, \boldsymbol{\lambda%
}\in\R^m \text{ or } \boldsymbol{\lambda}\in {\R^*}^m\right\}
\end{equation*}
and the induced group on the parameter is 
\begin{equation*}
\overline{\mathcal{G}}=\left\{\overline{g}:\boldsymbol{\theta}%
\in\Theta\mapsto \overline{g}(\boldsymbol{\theta}) = \text{diag}%
(\lambda_1,\ldots,\lambda_m) \boldsymbol{\theta};\, \boldsymbol{\lambda}\in\R%
^m \text{ or } \boldsymbol{\lambda}\in {\R^*}^m\right\}.
\end{equation*}
\end{example}

\vskip 0.5cm

In all the sequel, without loss of generality, we assume that the model $(%
\ref{the model})$ and the group of transformation $\mathcal{G}$ are such
that 
\begin{equation}  \label{assumption 1}
\mathbf{f}\left(g(\boldsymbol{x}),\overline{g}(\boldsymbol{\theta})\right) = 
\mathbf{f}\left(\boldsymbol{x},\boldsymbol{\theta}\right), \, \forall 
\boldsymbol{x}\in\R^m, \forall \boldsymbol{\theta}\in\Theta, \forall g\in 
\mathcal{G}. \\
\end{equation}
Note that this assumption implies the condition (\ref{invariance condition})
that the model $\mathcal{M}$ is invariant under $\mathcal{G}$. \newline

In all the following, when estimating $\boldsymbol{\theta}_T$ by an estimate 
$\widetilde{\boldsymbol{\theta}}_n := T(\boldsymbol{X}_1,\ldots,\boldsymbol{X%
}_n)$, we consider the quadratic loss function 
\begin{equation}  \label{L2 loss}
L_2\left(\widetilde{\boldsymbol{\theta}}_n,\boldsymbol{\theta}_T\right) := 
\mathbb{E}\left(\left\|\widetilde{\boldsymbol{\theta}}_n- \boldsymbol{\theta}%
_T\right\|^2\right)
\end{equation}
if the model is invariant with respect to additive group, and the loss
function is taken to be relative quadratic 
\begin{equation}  \label{Lr loss}
L_r\left(\widetilde{\boldsymbol{\theta}}_n,\boldsymbol{\theta}_T\right) :=%
\mathbb{E}\left(\left\|\frac{\widetilde{\boldsymbol{\theta}}_n}{\boldsymbol{%
\theta}_T}-\boldsymbol{1}\right \|^2\right)
\end{equation}
if the model is invariant with respect to the multiplicative group.\newline

\begin{definition}
(invariant loss under a group of transformations). \label{def loss
invariance} A loss function $L(\cdot,\cdot) : (\boldsymbol{\theta}_T, d)\in
\Theta\times D \mapsto L(\boldsymbol{\theta}_T,d)\in\R_+$, where $D$ denotes
the set of the parameter estimates (called decision space), is invariant
under a transformation $g\in\mathcal{G}$ iff for any estimate $d\in D$,
there exists a unique $d^{\prime }\in D$ such that $L(\boldsymbol{\theta}%
_T,d)=L(\overline{g}(\boldsymbol{\theta}_T),d^{\prime }).$ We denote then by 
$\widetilde{g}(\cdot)$ the bijection, from $D$ onto $D$, such that $%
d^{\prime }=\widetilde{g}(d)$. Hence, we have 
\begin{equation*}
L(\boldsymbol{\theta}_T,d) = L(\overline{g}(\boldsymbol{\theta}_T),%
\widetilde{g}(d)), \quad \forall d\in D.
\end{equation*}
\end{definition}

We denote by $\widetilde{\mathcal{G}}$ the induced group on the decision
space $D$.

\begin{definition}
Assume that the estimation problem $(\mathcal{M}, D, L)$ is invariant under
the group $\mathcal{G}$. Let $\overline{\mathcal{G}}$ and $\widetilde{%
\mathcal{G}}$ be, respectively, the induced groups on the parameter space $%
\Theta$ and the decision space $D$. An estimate $\widetilde{\boldsymbol{%
\theta}}_n:=T(\boldsymbol{X}_1,\ldots,\boldsymbol{X}_n)\in D$ is said to be
equivariant iff 
\begin{equation*}
T(g(\boldsymbol{X}_1),\ldots,g(\boldsymbol{X}_n)) = \widetilde{g}(T(%
\boldsymbol{X}_1,\ldots,\boldsymbol{X}_n)), \quad \forall g\in\mathcal{G}.
\end{equation*}
\end{definition}

We will see, under condition (\ref{assumption 1}), that the empirical
minimum divergence estimates, introduced in \cite{Bronia_Keziou2012}, are
equivariant for the above models, using results on the existence and
characterization of the distribution $P_{\boldsymbol{X}}$ ont the sets $%
\mathcal{M}_{\boldsymbol{\theta}}$. First, we recall the definition of $%
D_\varphi$-divergences and some of their properties. Let $\varphi$ be a
convex function from $\mathbb{R}$ onto $[0,+\infty]$ with $\varphi(1)=0$,
and such that its domain, $\text{dom}\varphi :=\left\{x\in\mathbb{R} \text{
such that } \varphi(x)<\infty\right\}=:(a,b),$ is an interval, with
endpoints satisfying $a<1<b$, which may be bounded or unbounded, open or
not. We assume that $\varphi$ is closed; the closedness of $\varphi$ means
that if $a$ or $b$ are finite then $\varphi(x)\to\varphi(a)$ when $x
\downarrow a$, and $\varphi(x)\to\varphi(b)$ when $x \uparrow b$. Note that,
this is equivalent to the fact that the level sets $\{x\in\mathbb{R};~
\varphi(x)\leq \alpha \}$, $\forall \alpha\in\mathbb{R}$, are closed in $%
\mathbb{R}$ endowed with the usual topology. For any s.f.m. $Q\in M$, the $%
D_\varphi$-divergence between $Q$ and a probability distribution $P$, when $%
Q $ is absolutely continuous with respect to (a.c.w.r.t) $P$, is defined
through 
\begin{equation}
D_{\varphi}(Q,P):=\int_{\mathbb{R}^m}\varphi\left(\frac{dQ}{dP}(\boldsymbol{x%
})\right)\, dP(\boldsymbol{x}),  \label{divRusch}
\end{equation}
where $\frac{dQ}{dP}(\cdot)$ is the Radon-Nikodym derivative of $Q$ w.r.t. $%
P $. When $Q$ is not a.c.w.r.t. $P$, we set $D_{\varphi}(Q,P):=+\infty$. For
any probability distribution $P$, the mapping $Q\in M\mapsto
D_{\varphi}(Q,P) $ is convex and takes nonnegative values. When $Q=P$ then $%
D_{\varphi}(Q,P)=0 $. Furthermore, if the function $x\mapsto\varphi(x)$ is
strictly convex on a neighborhood of $x=1$, then we have 
\begin{equation}
D_{\varphi}(Q,P)=0~\text{ if and only if }~Q=P.  \label{p.f.}
\end{equation}
\noindent All the above properties are presented in \cite{Csiszar1963}, \cite%
{Csiszar1967} and in Chapter 1 of \cite{Liese-Vajda1987}, for $D_\varphi-$%
divergences defined on the set of all probability distributions $M^{1}$.
When the $D_\varphi$-divergences are extended to $M$, then the same
arguments as developed on $M^{1}$ hold. When defined on $M^{1}$, the
Kullback-Leibler $(KL)$, modified Kullback-Leibler $(KL_{m})$, $\chi^{2}$,
modified $\chi^{2}$ $(\chi_{m}^{2})$, Hellinger $(H)$, and $L^{1}$
divergences are respectively associated to the convex functions $%
\varphi(x)=x\log x-x+1$, $\varphi(x)=-\log x+x-1$, $\varphi (x)=\frac{1}{2}{%
(x-1)}^{2}$, $\varphi(x)=\frac{1}{2}{(x-1)}^{2}/x $, $\varphi(x)=2{(\sqrt{x}%
-1)}^{2}$ and $\varphi(x)=\left\vert x-1\right\vert $. All these divergences
except the $L^{1}$ one, belong to the class of the so-called power
divergences introduced in \cite{Cressie-Read1984} (see also \cite%
{Liese-Vajda1987} and \cite{PardoLeandro2006}). They are defined through the
class of convex functions 
\begin{equation}
x\in\mathbb{R}_{+}^{\ast}\mapsto\varphi_{\gamma}(x):=\frac{x^{\gamma}-\gamma
x+\gamma-1}{\gamma(\gamma-1)}  \label{gamma convex functions}
\end{equation}
if $\gamma\in\mathbb{R}\setminus\left\{ 0,1\right\} $, $\varphi
_{0}(x):=-\log x+x-1$ and $\varphi_{1}(x):=x\log x-x+1$. So, the $KL-$%
divergence is associated to $\varphi_{1}$, the $KL_{m}$ to $\varphi_{0}$,
the $\chi^{2}$ to $\varphi_{2}$, the $\chi_{m}^{2}$ to $\varphi_{-1}$ and
the Hellinger distance to $\varphi_{1/2}$. We extend the definition of the
power divergences functions $Q\in M^{1}\mapsto D_{\varphi_{\gamma}}(Q,P)$
onto the whole set of signed finite measures $M$ as follows. When the
function $x\mapsto\varphi_{\gamma}(x)$ is not defined on $]-\infty,0[$ or
when $\varphi_{\gamma}$ is defined on $\mathbb{R}$ but is not convex (for
instance if $\gamma=3$), we extend the definition of $\varphi_{\gamma}$ as
follows 
\begin{equation}
x\in\mathbb{R}\mapsto\varphi_{\gamma}(x)\mathds{1}_{[0,+\infty[}(x)+(+\infty
)\mathds{1}_{]-\infty,0[}(x).  \label{gamma convex functions sur R}
\end{equation}
Note that for $\chi^2$-divergence, the corresponding $\varphi$ function $%
\varphi(x)=\frac{1}{2}(x-1)^2$ is convex and defined on whole $\mathbb{R}$.
In this paper, for technical considerations, we assume that the functions $%
\varphi$ are strictly convex on their domain $(a,b)$, twice continuously
differentiable on $]a,b[$, the interior of their domain. Hence, $%
\varphi^{\prime }(1)=0$, and for all $x\in ]a,b[$, $\varphi^{\prime \prime
}(x)>0$. Here, $\varphi^{\prime }$ and $\varphi^{\prime \prime }$ are used
to denote respectively the first and the second derivative functions of $%
\varphi$. Note that the above assumptions on $\varphi$ are not restrictive,
and that all the power functions $\varphi_{\gamma}$, see (\ref{gamma convex
functions sur R}), satisfy the above conditions, including all standard
divergences.

\section{Minimum empirical divergence estimates}

\noindent Let $\boldsymbol{X}_1,\ldots,\boldsymbol{X}_n$ denote an i.i.d.
sample of a random vector $\boldsymbol{X}\in\R^m$ with probability
distribution $P_{\mathbf{X}}$. Let $P_{n}(\cdot)$ be the associated
empirical measure, namely, 
\begin{equation*}
P_n(\cdot):=\frac{1}{n}\sum_{i=1}^n\delta_{\boldsymbol{X}_i}(\cdot),
\end{equation*}
where $\delta_{\boldsymbol{x}}(\cdot)$ denotes the Dirac measure at point $%
\boldsymbol{x}$, for all $\boldsymbol{x}$. For a given $\theta\in\Theta$,
the ``plug-in'' estimate of $D_{\varphi}(\mathcal{M}_{\theta},P_{\mathbf{X}%
}) $ is 
\begin{equation}  \label{plug-in estimate of phi}
\widetilde{D}_{\varphi}(\mathcal{M}_{\boldsymbol{\theta}},P_{\boldsymbol{X}%
}) := \inf_{Q\in\mathcal{M}_{\boldsymbol{\theta} }}D_{\varphi}(Q,P_n) =
\inf_{Q\in\mathcal{M}_{\boldsymbol{\theta}}}\int_{\R^m} \varphi\left(\frac {%
dQ}{dP_n}(\boldsymbol{x})\right)\, dP_n(\boldsymbol{x}).
\end{equation}
If the projection $Q^{(n)}_{\boldsymbol{\theta}}$ of $P_n$ on $\mathcal{M}_{%
\boldsymbol{\theta}}$ exists, then it is clear that $Q^{(n)}_{\boldsymbol{%
\theta}}$ is a s.f.m. (or possibly a probability distribution) a.c.w.r.t. $%
P_n$; this means that the support of $Q^{(n)}_{\boldsymbol{\theta}}$ must be
included in the set $\left\{\boldsymbol{X}_1,\ldots,\boldsymbol{X}_n\right\}$%
. So, define the set 
\begin{equation}  \label{m theta n}
\mathcal{M}_{\boldsymbol{\theta}}^{(n)}:=\left\{ Q\in M~ s.t. ~Q \text{
a.c.w.r.t. } P_{n},~\sum_{i=1}^{n}Q(\boldsymbol{X}_i)=1 \text{ and }
\sum_{i=1}^nQ(\boldsymbol{X}_i )\mathbf{f}(\boldsymbol{X}_i,\boldsymbol{%
\theta}) = \boldsymbol{0}\right\},
\end{equation}
which may be seen as a subset of $\R^n$. Then, the plug-in estimate (\ref%
{plug-in estimate of phi}) can be written as 
\begin{equation}  \label{estim de phiMtheta cas cont}
\widetilde{D}_{\varphi}(\mathcal{M}_{\boldsymbol{\theta}},P_{\boldsymbol{X}%
}) := \inf_{Q\in\mathcal{M}_{\boldsymbol{\theta}}}D_{\varphi}(Q,P_n) =
\inf_{Q\in\mathcal{M}_{\boldsymbol{\theta} }^{(n)}}\frac{1}{n}%
\sum_{i=1}^n\varphi\left( nQ(\boldsymbol{X}_i)\right).
\end{equation}
In the same way, 
\begin{equation*}
D_{\varphi}(\mathcal{M},P_{\boldsymbol{X}}):= \inf_{\boldsymbol{\theta}%
\in\Theta} \inf_{Q\in\mathcal{M}_{\boldsymbol{\theta}}}D_{\varphi}(Q,P_{%
\boldsymbol{X}})
\end{equation*}
can be estimated by 
\begin{equation}  \label{estim de phiM  cas cont}
\widetilde{D}_{\varphi}(\mathcal{M},P_{\boldsymbol{X}}):= \inf_{\boldsymbol{%
\theta}\in\Theta}\widetilde{D}_{\varphi}(\mathcal{M}_{\boldsymbol{\theta}%
},P_{\boldsymbol{X}}) = \inf_{\boldsymbol{\theta}\in\Theta} \inf_{Q\in%
\mathcal{M}_{\boldsymbol{\theta} }}D_{\varphi}(Q,P_n) =\inf_{\boldsymbol{%
\theta}\in\Theta}\inf _{Q\in\mathcal{M}_{\boldsymbol{\theta}}^{(n)}}\frac{1}{%
n}\sum_{i=1}^n\varphi\left( nQ(\boldsymbol{X}_i)\right).
\end{equation}
By uniqueness of $\arg\inf_{\boldsymbol{\theta}\in\Theta}D_{\varphi}(%
\mathcal{M}_{\boldsymbol{\theta}},P_{\boldsymbol{X}})$ and since the infimum
is reached in $\boldsymbol{\theta}=\boldsymbol{\theta}_T$, we estimate $%
\boldsymbol{\theta}_T$ through 
\begin{equation}  \label{theta estimate}
\widetilde{\boldsymbol{\theta}}_{\varphi} :=\arg\inf_{\boldsymbol{\theta}%
\in\Theta}\inf_{Q\in\mathcal{M}_{\boldsymbol{\theta}}^{(n)}} \frac{1}{n}%
\sum_{i=1}^n\varphi\left(nQ(\boldsymbol{X}_i)\right).
\end{equation}
The expression of the estimate $\widetilde{D}(\mathcal{M}_{\boldsymbol{\theta%
}},P_{\boldsymbol{X}})$, given in (\ref{estim de phiMtheta cas cont}), is
the solution of a convex optimization problem under convex constrained
subset in $\mathbb{R}^n$. In order to transform this problem to an
unconstrained one, we will make use of the Fenchel-Lengendre transform,
denoted $\varphi^*(\cdot)$, of the convex function $\varphi(\cdot)$, as well
as some other duality arguments. It is defined by 
\begin{equation}
\varphi^* : t\in \mathbb{R}\mapsto \varphi^*(t):=\sup_{x\in\mathbb{R}%
}\left\{tx-\varphi(x)\right\}.
\end{equation}

\noindent For convenience, we recall some properties of the convex conjugate 
$\varphi^*$ of $\varphi$. For the proofs we can refer to Section 26 in \cite%
{Rockafellar1970}. Theses properties will be used to determine the convex
conjugates $\varphi^*$ of some standard divergence functions $\varphi$; see
Table \ref{table convex conjugates} below. The function $\varphi^*$ in turn
is convex and closed, its domain is an interval $(a^*,b^*)$ with endpoints 
\begin{equation}  \label{domain de varphi*}
a^* = \lim_{x\to -\infty}\frac{\varphi(x)}{x},\quad b^* = \lim_{x \to
+\infty}\frac{\varphi(x)}{x}
\end{equation}
satisfying $a^*<0<b^*$ with $\varphi^*(0)=0$. Note that the interval 
\begin{equation*}
(\widetilde{a}, \widetilde{b}):=\left(\lim_{x\downarrow a}\frac{\varphi(x)}{x%
},\lim_{x\uparrow b}\frac{\varphi(x)}{x}\right)
\end{equation*}
can be different from $(a^*,b^*)$, the real domain of $\varphi^*$ given by (%
\ref{domain de varphi*}). This holds when $a$ or $b$ is finite and $%
\varphi^{\prime }(a)$ or $\varphi^{\prime }(b)$ is finite, respectively. For
example, for the convex function 
\begin{equation*}
\varphi(x)=\frac{1}{2}(x-1)^2\mathds{1}_{\mathbb{R}_+}(x)+(+\infty)\mathds{1}%
_{]0,+\infty[}(x),
\end{equation*}
we have $\text{dom}\varphi =[0,+\infty[$ and $\varphi^{\prime }(0)=-1$, and
we can see that the domain of the corresponding $\varphi^*$-function is $%
(a^*,b^*)=(-\infty, +\infty)$ which is different from $(\widetilde{a}, 
\widetilde{b}) :=\left(\lim_{x\downarrow 0}\frac{\varphi(x)}{x},\lim_{x
\uparrow +\infty}\frac{\varphi(x)}{x}\right)=(-1,+\infty).$ The two
intervals $(a^*,b^*)$ and $(\widetilde{a}, \widetilde{b})$ coincide if the
function $\varphi$ is ``essentially smooth'', i.e., differentiable with 
\begin{equation}
\begin{array}{ccccc}
\lim_{t\downarrow a}{\varphi}^{\prime}(t) & = & -\infty & \text{ if } & a 
\text{ is finite}, \\ 
\lim_{t\uparrow b}{\varphi}^{\prime}(t) & = & +\infty & \text{ if } & b 
\text{ is finite}.%
\end{array}%
\end{equation}
The strict convexity of $\varphi$ on its domain $(a,b)$ is equivalent to the
condition that its conjugate $\varphi^*$ is essentially smooth, i.e.,
differentiable with 
\begin{equation}
\begin{array}{ccccc}
\lim_{t\downarrow a^*}{\varphi^*}^{\prime}(t) & = & -\infty & \text{ if } & 
a^* \text{ is finite}, \\ 
\lim_{t\uparrow b^*}{\varphi^*}^{\prime}(t) & = & +\infty & \text{ if } & 
b^* \text{ is finite}.%
\end{array}%
\end{equation}
Conversely, $\varphi$ is essentially smooth on its domain $(a,b)$ if and
only if $\varphi^*$ is strictly convex on its domain $(a^*,b^*)$. \newline

\noindent In all the sequel, we assume additionally that $\varphi$ is
essentially smooth. Hence, $\varphi^*$ is strictly convex on its domain $%
(a^*,b^*)$, and it holds that 
\begin{equation*}
a^* =\lim_{x\to -\infty}\frac{\varphi(x)}{x}=\lim_{x\downarrow a}\frac{%
\varphi(x)}{x} =\lim_{x\downarrow a}\varphi^{\prime * }=\lim_{x\to + \infty}%
\frac{\varphi(x)}{x}=\lim_{x\uparrow b}\frac{\varphi(x)}{x}= \lim_{x\uparrow
b}\varphi^{\prime }(x),
\end{equation*}
and 
\begin{equation}  \label{forme explicite de psi}
\varphi^*(t)=t{\varphi^{\prime }}^{-1}(t)-\varphi\left( {\varphi^{\prime }}%
^{-1}(t)\right),\quad \text{for all } t\in ]a^*,b^*[,
\end{equation}
where ${\varphi^{\prime }}^{-1}$ denotes the inverse of the derivative
function $\varphi^{\prime }$ of $\varphi$. It holds also that $\varphi^*$ is
twice continuously differentiable on $]a^*,b^*[$ with 
\begin{equation}  \label{derivee de psi}
{\varphi^*}^\prime(t)={\varphi^{\prime }}^{-1}(t) \quad \text{and}\quad {%
\varphi^*}^{\prime\prime}(t)=\frac{1}{\varphi^{\prime \prime }\left({%
\varphi^{\prime }}^{-1}(t)\right)}.
\end{equation}
In particular, ${\varphi^*}^\prime(0)=1$ and ${\varphi^*}^{\prime%
\prime}(0)=1 $. Obviously, since $\varphi$ is assumed to be closed, we have 
\begin{equation*}
\varphi(a)=\lim_{x\downarrow a}\varphi(x)\quad \text{ and }\quad
\varphi(b)=\lim_{x\uparrow b}\varphi(x),
\end{equation*}
which may be finite or infinite. Hence, by closedness of $\varphi^*$,
likewise we have 
\begin{equation*}
\varphi^*(a^*)=\lim_{t\downarrow a^*}\varphi^*(x)\quad \text{ and }\quad
\varphi^*(b^*)=\lim_{t\uparrow b^*}\varphi^*(t).
\end{equation*}
Finally, the first and second derivatives of $\varphi$ in $a$ and $b$ are
defined to be the limits of $\varphi^{\prime }(x)$ and $\varphi^{\prime
\prime }(x)$ when $x\downarrow a$ and when $x\uparrow b$. The first and
second derivatives of $\varphi^*$ in $a^*$ and $b^*$ are defined in a
similar way. In Table \ref{table convex conjugates}, using the above
properties, we give the convex conjugates $\varphi^*$ of some standard
divergence functions $\varphi$, associated to standard divergences. We
determine also their domains, respectively, $(a,b)$ and $(a^*,b^*)$.

\begin{table}[h]
\caption{Convex conjugates $\protect\varphi^*$ of some standard divergence
functions $\protect\varphi$.}
\label{table convex conjugates}%
\begin{tabular}{|l||l|l||l|l|}
\hline
$D_\varphi$ & $\text{dom}\varphi =:(a,b) $ & $\varphi$ & $\text{dom}%
\varphi^* =:(a^*,b^*)$ & $\varphi^*$ \\ \hline\hline
$D_{KL_m}$ & $]0,+\infty[$ & $\varphi(x):=-\log x +x -1$ & $]-\infty,1[$ & $%
\varphi^*(t)= - \log(1-t)$ \\ \hline
$D_{KL}$ & $[0,+\infty[$ & $\varphi(x):=x\log x -x +1$ & $\mathbb{R}$ & $%
\varphi^*(t)=e^t-1$ \\ \hline
$D_{\chi^2_m}$ & $]0,+\infty[$ & $\varphi(x):= \frac{1}{2}\frac{%
\left(x-1\right)^2}{x}$ & $\left]-\infty,\frac{1}{2}\right]$ & $%
\varphi^*(t)=1-\sqrt{1-2t}$ \\ \hline
$D_{\chi^2}$ & $\mathbb{R}$ & $\varphi(x):= \frac{1}{2}\left(x-1\right)^2$ & 
$\mathbb{R}$ & $\varphi^*(t)=\frac{1}{2}t^2+t$ \\ \hline
$D_H$ & $[0,+\infty[$ & $\varphi(x):=2(\sqrt{x}-1)^2$ & $]-\infty, 2[$ & $%
\varphi^*(t)=\frac{2t}{2-t}$ \\ \hline
\end{tabular}%
\end{table}

\noindent Using some duality arguments, see \cite{Bronia_Keziou2012}, we can
show that, for any $\boldsymbol{\theta}\in\Theta$, if there exists $Q_0$ in $%
\mathcal{M}_{\boldsymbol{\theta}}^{(n)}$ such that 
\begin{equation}  \label{cond Owen}
a< nQ_0(\boldsymbol{X}_i)<b, \quad \text{for all} \quad i=1,\ldots,n,
\end{equation}
then 
\begin{equation}  \label{egalite duale 1}
\widetilde{D}_{\varphi}(\mathcal{M}_{\boldsymbol{\theta}},P_{\boldsymbol{X}%
}):=\inf_{Q\in\mathcal{M}_{\boldsymbol{\theta}}}D_{\varphi}(Q,P_n) =
\sup_{(t_0,t_1,\ldots,t_\ell)^\top\in \R^{1+\ell}}\left\{t_0 - \frac{1}{n}%
\sum_{i=1}^n\varphi^*( t_0 +\sum_{j=1}^\ell t_jf_j(\boldsymbol{X}_i,%
\boldsymbol{\theta})) \right\}
\end{equation}
with dual attainment. Conversely, if there exists some dual optimal solution 
\begin{equation*}
\widetilde{\underline{\boldsymbol{t}}}(\boldsymbol{\theta}):=\left(%
\widetilde{t}_0(\boldsymbol{\theta}),\widetilde{t}_1(\boldsymbol{\theta}%
),\ldots,\widetilde{t}_\ell(\boldsymbol{\theta})\right)^{\top}\in\mathbb{R}%
^{1+\ell}
\end{equation*}
such that 
\begin{equation}  \label{cond Owen 2}
a^*< \widetilde{t_0}(\boldsymbol{\theta})+\sum_{j=1}^\ell \widetilde{t_j}(%
\boldsymbol{\theta}) f_j(\boldsymbol{X}_i,\boldsymbol{\theta})<b^*, \quad 
\text{for all}\quad i=1,\ldots,n,
\end{equation}
then the equality (\ref{egalite duale 1}) holds, and the unique optimal
solution of the primal problem 
\begin{equation*}
\inf_{Q\in\mathcal{M}_{\boldsymbol{\theta}}^{(n)}}D_{\varphi}(Q,P_n),
\end{equation*}
namely, the projection of $P_n$ on $\mathcal{M}^{(n)}_{\boldsymbol{\theta}}$%
, is given by 
\begin{equation}  \label{estimateur du projete}
Q^{(n)}_{\boldsymbol{\theta}}(\boldsymbol{X}_i)=\frac{1}{n} {\varphi^{\prime
}}^{-1}\left(\widetilde{t_0}(\boldsymbol{\theta})+\sum_{j=1}^{\ell}%
\widetilde{t_j}(\boldsymbol{\theta}) f_j(\boldsymbol{X}_i,\boldsymbol{\theta}%
)\right),\quad i=1,\ldots,n,
\end{equation}
where $\widetilde{\underline{\boldsymbol{t}}}(\boldsymbol{\theta}):=\left(%
\widetilde{t}_0(\boldsymbol{\theta}),\widetilde{t}_1(\boldsymbol{\theta}%
),\ldots,\widetilde{t}_\ell(\boldsymbol{\theta})\right)^{\top}\in \mathbb{R}%
^{1+\ell}$ is solution of the system of equations 
\begin{equation}  \label{equation-system-n}
\left\{%
\begin{array}{lll}
1-\frac{1}{n}\sum_{i=1}^{n} {\varphi^{\prime }}^{-1}(t_0 +\sum_{j=1}^\ell
t_j f_j (\boldsymbol{X}_i,\boldsymbol{\theta})) & = & 0, \\ 
-\frac{1}{n}\sum_{i=1}^n f_{j}(\boldsymbol{X}_i,\boldsymbol{\theta}) {%
\varphi^{\prime }}^{-1}( t_0 +\sum_{j=1}^{\ell} t_j f_j(\boldsymbol{X}_i,%
\boldsymbol{\theta})) & = & 0, \quad j=1,\ldots,\ell. \\ 
&  & 
\end{array}
\right.
\end{equation}
In view of the last results, using the notations 
\begin{equation*}
\underline{\boldsymbol{t}} :=(t_0,t_1,\ldots,t_\ell)^\top \in\R^{1+\ell},
\quad \underline{\mathbf{f}} (\cdot,\boldsymbol{\theta}) := \left(\mathds{1}%
_{\R^m}(\cdot), f_1(\cdot, \boldsymbol{\theta}), \ldots, f_\ell(\cdot,%
\boldsymbol{\theta})\right)^\top
\end{equation*}
and 
\begin{equation*}
\widetilde{\underline{\boldsymbol{t}}}(\boldsymbol{\theta}) :=\arg\sup_{%
\underline{\boldsymbol{t}}\in\mathbb{R}^{1+\ell}} \left\{ t_{0}-\frac{1}{n}%
\sum_{i=1}^{n}\varphi^*\left( \underline{\boldsymbol{t}}^\top \underline{%
\mathbf{f}}(X_{i},\boldsymbol{\theta})\right)\right\}, \quad \forall 
\boldsymbol{\theta}\in\Theta,
\end{equation*}
we obtain the following equivalent expressions to the estimates $\widetilde{D%
}_{\varphi}(\mathcal{M}_{\boldsymbol{\theta}},P_{\mathbf{X}})$, $\widetilde{D%
}_{\varphi}(\mathcal{M},P_{\boldsymbol{X}})$ and $\widetilde{\boldsymbol{%
\theta}}$, see (\ref{plug-in estimate of phi}), (\ref{estim de phiM cas cont}%
) and (\ref{theta estimate}), 
\begin{eqnarray}  \label{estimate m theta 1}
\widetilde{D}_{\varphi}(\mathcal{M}_{\theta},P_{\mathbf{X}}) & = & \sup_{%
\underline{\boldsymbol{t}}\in\mathbb{R}^{1+\ell}}\left\{ t_{0}-\frac{1}{n}%
\sum_{i=1}^{n}\varphi^*\left( \underline{\boldsymbol{t}}^\top \underline{%
\mathbf{f}}(X_{i},\boldsymbol{\theta})\right)\right\} \\
& = & \widetilde{t}_0(\boldsymbol{\theta})-\frac{1}{n}\sum_{i=1}^{n}%
\varphi^*\left( \widetilde{\underline{\boldsymbol{t}}}(\boldsymbol{\theta}%
)^\top \underline{\mathbf{f}}(X_{i},\boldsymbol{\theta})\right),
\label{estimate m theta 2}
\end{eqnarray}
\begin{eqnarray}  \label{estimate m}
\widetilde{D}_{\varphi}(\mathcal{M},P_{\boldsymbol{X}}):=\inf_{\boldsymbol{%
\theta}\in\Theta} \widetilde{D}_{\varphi}(\mathcal{M}_{\boldsymbol{\theta}%
},P_{\boldsymbol{X}}) & = & \inf_{\boldsymbol{\theta}\in\Theta}\sup_{%
\underline{\boldsymbol{t}}\in \mathbb{R}^{1+\ell}}\left\{ t_{0}-\frac{1}{n}%
\sum_{i=1}^{n}\varphi^*\left( \underline{\boldsymbol{t}}^\top \underline{%
\mathbf{f}}(\boldsymbol{X}_i,\boldsymbol{\theta})\right)\right\} \\
& = & \inf_{\boldsymbol{\theta}\in\Theta} \left\{ \widetilde{t}_0(%
\boldsymbol{\theta})-\frac{1}{n}\sum_{i=1}^{n}\varphi^*\left( \widetilde{%
\underline{\boldsymbol{t}}}(\boldsymbol{\theta})^\top \underline{\mathbf{f}}%
(X_{i},\boldsymbol{\theta})\right)\right\},
\end{eqnarray}
and 
\begin{equation}  \label{estimate theta}
\widetilde{\boldsymbol{\theta}}_\varphi = \arg\inf_{\boldsymbol{\theta}%
\in\Theta} \left\{ \widetilde{t}_0(\boldsymbol{\theta})-\frac{1}{n}%
\sum_{i=1}^{n}\varphi^*\left( \widetilde{\underline{\boldsymbol{t}}}(%
\boldsymbol{\theta})^\top \underline{\mathbf{f}}(X_{i},\boldsymbol{\theta}%
)\right)\right\}. \\
\end{equation}

\begin{remark}
The empirical likelihood estimate $\widetilde{\theta}_{EL}$ is obtained for
the particular choice of the modified Kullback-Leibler divergence $KL_m$,
namely, when $\varphi(x) = \varphi_0(x) =-\log x + x -1$. Moreover,
straightforward computation shows that $\widetilde{t}_0(\boldsymbol{\theta}%
)=0$, $\forall \boldsymbol{\theta}\in\Theta$. Therefore, $t_0$ can be
omitted, and the above expression can be simplified to 
\begin{equation*}
\widetilde{D}_{KL_m}(\mathcal{M},P_{\boldsymbol{X}}) = \inf_{\boldsymbol{%
\theta}\in\Theta}\sup_{\boldsymbol{t}\in\R^\ell} \frac{1}{n}\sum_{i=1}^n
\log\left(1-\boldsymbol{t}^\top\mathbf{f} (\boldsymbol{X}_i,\boldsymbol{%
\theta})\right)
\end{equation*}
and 
\begin{equation*}
\widetilde{\boldsymbol{\theta}}_{EL}=\widetilde{\boldsymbol{\theta}}%
_{\varphi_0} = \arg\inf_{\boldsymbol{\theta}\in\Theta} \sup_{\boldsymbol{t}%
\in\R^\ell} \frac{1}{n}\sum_{i=1}^n \log\left(1-\boldsymbol{t}^\top\mathbf{f}
(\boldsymbol{X}_i,\boldsymbol{\theta})\right).
\end{equation*}%
\newline
\end{remark}

We will show that for any divergence $D_\varphi$, the estimate $\widetilde{%
\boldsymbol{\theta}}_\varphi := \widetilde{\boldsymbol{\theta}}_ \varphi(%
\boldsymbol{X}_1,\ldots,\boldsymbol{X}_n)$ is invariant with respect to $L_2$
loss for the additive group, and invariant with respect to $L_r$ loss for
the multiplicative group. First, we expose the asymptotic counterpart of the
estimates (\ref{estimate m theta 1}), (\ref{estimate m}) and (\ref{estimate
theta}). In particular, we give results about existence and characterization
of the projection of $P_X$ on the model $\mathcal{M}$. The characterization
of the projection will be of great importance in computing the minimum risk
equivariant estimate. We have; see Theorem 1 in \cite{Bronia_Kez2006_STUDIA}:

\begin{proposition}
\label{proposition 2} \label{Prop caract theo} Let $\boldsymbol{\theta}$ be
a given value in $\Theta$. Assume that $\int_{\R^m} |f_j(\boldsymbol{x},%
\boldsymbol{\theta})|\,dP_{\boldsymbol{X}}(\boldsymbol{x})<\infty$ for all $%
j=1,\ldots,\ell$, and that there exists $Q_0$ in $\mathcal{M}_{\boldsymbol{%
\theta}}$ such that $D_\varphi(Q_0,P_{\boldsymbol{X}})<\infty$ and\footnote{%
The strict inequalities mean that $P_{\boldsymbol{X}}\left(\left\{%
\boldsymbol{x}\in\R^m;\, \frac{dQ_0}{dP_{\boldsymbol{X}}}(\boldsymbol{x}%
)\leq a\right\}\right)=P_{\boldsymbol{X}}\left(\left\{\boldsymbol{x}\in\R%
^m;\,\frac{dQ_0}{dP_{\boldsymbol{X}}}(\boldsymbol{x})\geq
b\right\}\right)=0. $} 
\begin{equation}  \label{condition de qualification}
a<\inf_{\boldsymbol{x}\in\R^m}\frac{dQ_{0}}{dP_{\boldsymbol{X}}}(\boldsymbol{%
x})\leq\sup_{\boldsymbol{x}\in\R^m}\frac{dQ_0}{dP_{\boldsymbol{X}}}(%
\boldsymbol{x})<b, \quad P_{\boldsymbol{X}}-a.s.
\end{equation}
Then, we have 
\begin{equation}  \label{egalite duale 2}
\inf_{Q\in\mathcal{M}_{\boldsymbol{\theta}}}D_{\varphi}(Q,P_{\boldsymbol{X}%
})=\sup_{\underline{\boldsymbol{t}} \in\mathbb{R}^{1+\ell} }\left\{
t_0-\int_{\mathbb{R}^m} \varphi^*( \underline{\boldsymbol{t}}^\top 
\underline{\mathbf{f}} (\boldsymbol{x},\boldsymbol{\theta}))\,dP_{%
\boldsymbol{X}}(\boldsymbol{x})\right\}
\end{equation}
with dual attainment. Conversely, if there exists a dual optimal solution 
\begin{equation*}
\underline{\boldsymbol{t}}(\boldsymbol{\theta}) := (t_0(\boldsymbol{\theta}%
),t_1(\boldsymbol{\theta}),\ldots,t_\ell(\boldsymbol{\theta}))^\top\in\R%
^{1+\ell}
\end{equation*}
belonging to the interior (in $\R^{1+\ell}$) of the set 
\begin{equation}  \label{condition de qualification 2}
\left\{\underline{\boldsymbol{t}}\in\R^{1+\ell}~\text{such that}~ \int_{\R%
^m} |\varphi^*( \underline{\boldsymbol{t}}^\top\underline{\mathbf{f}}(%
\boldsymbol{x},\boldsymbol{\theta}))|\,dP_{\boldsymbol{X}}(\boldsymbol{x})
<\infty \right\},
\end{equation}
then the dual equality (\ref{egalite duale 2}) holds, and the unique optimal
solution $Q^*_{\boldsymbol{\theta}}$ of the primal problem $\inf_{Q\in%
\mathcal{M}_{\boldsymbol{\theta}}}D_{\varphi}(Q,P_{\boldsymbol{X}})$,
namely, the projection of $P_{\boldsymbol{X}}$ on $\mathcal{M}_{\boldsymbol{%
\theta}}$, is given by 
\begin{equation*}
\frac{dQ^*_{\boldsymbol{\theta}}}{dP_{\boldsymbol{X}}}(\boldsymbol{x})={%
\varphi^{\prime }}^{-1}\left(t_0(\boldsymbol{\theta})+\sum_{j=1}^\ell t_j(%
\boldsymbol{\theta})f_j(\boldsymbol{x},\boldsymbol{\theta})\right),
\end{equation*}
where $\underline{\boldsymbol{t}} (\boldsymbol{\theta}) := \left(t_0(%
\boldsymbol{\theta}),t_1(\boldsymbol{\theta}),\ldots,t_\ell(\boldsymbol{%
\theta})\right)^\top\in\R^{1+\ell}$ is the solution of the system of
equations 
\begin{equation}  \label{systeme}
\left\{%
\begin{array}{lll}
1-\int {\varphi^{\prime }}^{-1}( t_0+\sum_{j=1}^{\ell}t_j f_j(\boldsymbol{x},%
\boldsymbol{\theta}))\, dP_{\boldsymbol{X}}(\boldsymbol{x}) & = & 0 \\ 
- \int f_j(\boldsymbol{x},\boldsymbol{\theta}){\varphi^{\prime }}%
^{-1}(t_0+\sum_{j=1}^{\ell}t_j f_j(\boldsymbol{x},\boldsymbol{\theta}))\,
dP_{\boldsymbol{X}}(\boldsymbol{x}) & = & 0,\quad j =1,\ldots,\ell. \\ 
&  & 
\end{array}
\right.
\end{equation}
Furthermore, the solution $\underline{\boldsymbol{t}}(\boldsymbol{\theta})$
is unique if the functions $\mathds{1}_{\R^m}(\cdot), f_1(\cdot,\boldsymbol{%
\theta}),\ldots,f_\ell(\cdot,\boldsymbol{\theta})$ are linearly independent
in the sense that $P_{\boldsymbol{X}}\left(\left\{\boldsymbol{x}\in\R^m; \,
\, t_0+\sum_{j=1}^\ell t_jf_j(\boldsymbol{x},\boldsymbol{\theta})\neq
0\right\}\right)>0$ for all $\underline{\boldsymbol{t}} :=
(t_0,t_1,\ldots,t_\ell)^\top\in\R^{1+\ell}$ with $(t_0,t_1,\ldots,t_\ell)^%
\top\neq \boldsymbol{0}.$\newline
\end{proposition}

\begin{remark}
By minimizing $D_\varphi(Q,P_{\boldsymbol{X}})$, upon $Q\in\mathcal{M}_{%
\boldsymbol{\theta}}$, $\boldsymbol{\theta}\in\Theta$, we obtain the
semiparametric model of densities 
\begin{equation*}
p(\boldsymbol{x},\boldsymbol{\theta}) := {\varphi^{\prime }}^{-1}\left(t_0(%
\boldsymbol{\theta})+\sum_{j=1}^\ell t_j(\boldsymbol{\theta})f_j (%
\boldsymbol{x},\boldsymbol{\theta})\right) p_{\boldsymbol{X}}(\boldsymbol{x}%
);\, \boldsymbol{\theta}\in\Theta,
\end{equation*}
where $\underline{\boldsymbol{t}} (\boldsymbol{\theta}) := \left(t_0(%
\boldsymbol{\theta}),t_1(\boldsymbol{\theta}),\ldots,t_\ell(\boldsymbol{%
\theta})\right)^\top\in\R^{1+\ell}$ is the solution of the system of
equations (\ref{systeme}). For the particular case of the $KL$-divergence,
namely, when $\varphi(x)=\varphi_1(x)=x\log x -x+1$, $t_0(\boldsymbol{\theta}%
)$ can be explicitly computed, and the obtained model is the semiparametric
exponential family of probability densities 
\begin{eqnarray}  \label{exponential family}
p(\boldsymbol{x},\boldsymbol{\theta}) & := & \frac{\exp\left\{\sum_{j=1}^%
\ell t_j(\boldsymbol{\theta})f_j (\boldsymbol{x},\boldsymbol{\theta}%
)\right\}p_{\boldsymbol{X}}(\boldsymbol{x})}{\int_{\R^m}\exp\left\{%
\sum_{j=1}^\ell t_j(\boldsymbol{\theta})f_j (\boldsymbol{x}, \boldsymbol{%
\theta})\right\}\,dP_{\boldsymbol{X}}(\boldsymbol{x})}  \notag \\
& := & \frac{\exp\left\{\boldsymbol{t}(\boldsymbol{\theta})^\top \mathbf{f}(%
\boldsymbol{x},\boldsymbol{\theta})\right\} p_{\boldsymbol{X}}(\boldsymbol{x}%
)}{\int_{\R^m}\exp\left\{\boldsymbol{t}(\boldsymbol{\theta})^\top \mathbf{f}(%
\boldsymbol{x},\boldsymbol{\theta})\right\}\,dP_{\boldsymbol{X}}(\boldsymbol{%
x})} ;\, \boldsymbol{\theta}\in\Theta,
\end{eqnarray}
where, for all $\boldsymbol{\theta}\in\Theta$, $\boldsymbol{t}(\boldsymbol{%
\theta}) :=\left(t_1(\boldsymbol{\theta}),\ldots,t_\ell(\boldsymbol{\theta}%
)\right)^\top\in\R^\ell$ is the solution in $\boldsymbol{t}\in\R^\ell$ of
the system of equations 
\begin{equation}  \label{exponential family system}
\int_{\R^m}f_j(\boldsymbol{x},\boldsymbol{\theta}) \exp\left\{\boldsymbol{t}%
^\top \mathbf{f}(\boldsymbol{x},\boldsymbol{\theta})\right\}\, dP_{%
\boldsymbol{X}}(\boldsymbol{x})=0,\, j=1,\ldots,\ell.
\end{equation}
or equivalently 
\begin{equation*}
\int_{\R^m}\mathbf{f}(\boldsymbol{x},\boldsymbol{\theta}) \exp\left\{%
\boldsymbol{t}^\top \mathbf{f}(\boldsymbol{x},\boldsymbol{\theta})\right\}\,
dP_{\boldsymbol{X}}(\boldsymbol{x})=\mathbf{0}.
\end{equation*}
\end{remark}

\begin{proposition}
Assume that condition (\ref{assumption 1}) holds. Then, the minimum
empirical $\phi$-divergence estimates $(\ref{estimate theta})$ are
equivariant

\begin{enumerate}
\item[-] to the additive group of transformations with respect to the $L_2$
loss;

\item[-] to the multiplicative group of transformations with respect to the $%
L_r$ loss.
\end{enumerate}

Moreover, in both cases, the induced group of transformations $\widetilde{%
\mathcal{G}}$ on the space of estimates is equal to $\overline{\mathcal{G}}$%
, the group of transformations on the parameter space $\Theta$, in the sense
that 
\begin{equation*}
\widetilde{\boldsymbol{\theta}}_\varphi(g(\boldsymbol{X}_1),\ldots,g(%
\boldsymbol{X}_n)) = \overline{g}(\widetilde{\boldsymbol{\theta}}_\varphi(%
\boldsymbol{X}_1,\ldots,\boldsymbol{X}_n)), \, \forall g\in\mathcal{G}.
\end{equation*}
\end{proposition}

\begin{corollary}
For any estimate $\widetilde{\boldsymbol{\theta}}_\varphi$, the
corresponding loss function $\boldsymbol{\theta}_T\in\Theta \mapsto L(%
\boldsymbol{\theta}_T, \widetilde{\boldsymbol{\theta}}_\varphi)$ is constant.
\end{corollary}

In view of the above corollary, for the additive group, in order to obtain
the uniform minimum risk estimate, we can compute the risk $L_2(\boldsymbol{0%
},\widetilde{\boldsymbol{\theta}}_\varphi)$ of any estimate $\widetilde{%
\boldsymbol{\theta}}_\varphi$ under the particular value $\boldsymbol{\theta}%
_T=\boldsymbol{0}$, and then select the estimate that minimizes the risk.
Likewise, if a multiplicative group is considered, to obtain the uniform
minimum risk estimate, we can compute the risk $L_r(\mathbf{1},\widetilde{%
\boldsymbol{\theta}}_\varphi)$ of any estimate $\widetilde{\boldsymbol{\theta%
}}_\varphi$ under the particular value $\boldsymbol{\theta}_T=\mathbf{1}$,
and then select the estimate that minimizes the risk. To do this, we will
first characterize the equivariant estimates.

\begin{definition}
A functional $U : (\boldsymbol{x}_1,\ldots, \boldsymbol{x}_n)\in \R^{mn}
\mapsto U_{(\boldsymbol{x}_1,\ldots,\boldsymbol{x}_n)}(\cdot)\in \widetilde{%
\mathcal{G}}$ is ``invariant'' iff 
\begin{equation*}
U_{(g(\boldsymbol{x}_1),\ldots,g(\boldsymbol{x}_n))}(\cdot) = U_{(%
\boldsymbol{x}_1,\ldots,\boldsymbol{x}_n)}(\cdot), \, \forall g\in\mathcal{G}%
, \forall (\boldsymbol{x}_1,\ldots,\boldsymbol{x}_n) \in\R^{mn}.
\end{equation*}
\end{definition}

\begin{definition}
A functional $\,\mathcal{U} : (\boldsymbol{x}_1,\ldots,\boldsymbol{x}_n)\in %
\R^{mn} \mapsto \mathcal{U}_{(\boldsymbol{x}_1,\ldots,\boldsymbol{x}%
_n)}(\cdot)\in \widetilde{\mathcal{G}}$ is a ``maximal invariant'' iff it is
invariant and satisfies $\forall (\boldsymbol{x}_1,\ldots,\boldsymbol{x}_n),
\forall (\boldsymbol{y}_1,\ldots,\boldsymbol{y}_n)\in\R^{mn},$ 
\begin{equation*}
\mathcal{U}_{(\boldsymbol{x}_1,\ldots,\boldsymbol{x}_n)}(\cdot)= \mathcal{U}%
_{(\boldsymbol{y}_1,\ldots,\boldsymbol{y}_n)}(\cdot) \Rightarrow (%
\boldsymbol{y}_1,\ldots,\boldsymbol{y}_n) = (g(\boldsymbol{x}_1),\ldots,g(%
\boldsymbol{x}_n)),\, \text{ for some }g\in\mathcal{G}.
\end{equation*}
\end{definition}

\begin{remark}
For the additive group, we have that a functional $\mathcal{U}_{(\boldsymbol{%
x}_1,\ldots,\boldsymbol{x}_n)}(\cdot):=f_{(\boldsymbol{0},\boldsymbol{x}_2-%
\boldsymbol{x}_1,\ldots,\boldsymbol{x}_n-\boldsymbol{x}_1)}(\cdot)$, a
function of $(\boldsymbol{0},\boldsymbol{x}_2-\boldsymbol{x}_1,\ldots,%
\boldsymbol{x}_n-\boldsymbol{x}_1)$, is maximal invariant. Likewise, for the
multiplicative groupe, a functional $\mathcal{U}_{(\boldsymbol{x}_1,\ldots,%
\boldsymbol{x}_n)}(\cdot):=f_{\left(\frac{\boldsymbol{x}_1}{|\boldsymbol{x}%
_1|},\frac{\boldsymbol{x}_2} {\boldsymbol{x}_1},\ldots,\frac{\boldsymbol{x}_n%
}{\boldsymbol{x}_1}\right)}(\cdot)$, a function of $\left(\frac{\boldsymbol{x%
}_1}{|\boldsymbol{x}_1|},\frac{\boldsymbol{x}_2} {\boldsymbol{x}_1},\ldots,%
\frac{\boldsymbol{x}_n}{\boldsymbol{x}_1}\right)$, is maximal invariant.
\end{remark}

\begin{proposition}
Assume that the estimation problem $(\mathcal{M}, D, L)$ is invariant under
the group $\mathcal{G}$. Let $\overline{\mathcal{G}}$ and $\widetilde{%
\mathcal{G}}$ be, respectively, the induced groups on the parameter space $%
\Theta$ and the decision space $D$. Let $\widehat{\boldsymbol{\theta}}_0(%
\boldsymbol{X}_1,\ldots,\boldsymbol{X}_n)$ be any equivariant estimate.
Then, an estimate $\widetilde{\boldsymbol{\theta}}(\boldsymbol{X}_1,\ldots,%
\boldsymbol{X}_n)$ is equivariant iff 
\begin{equation*}
\widetilde{\boldsymbol{\theta}}(\boldsymbol{X}_1,\ldots,\boldsymbol{X}_n) =
U_{(\boldsymbol{X}_1,\ldots,\boldsymbol{X}_n)}\left(\widetilde{\boldsymbol{%
\theta}}_0(\boldsymbol{X}_1,\ldots,\boldsymbol{X}_n)\right),
\end{equation*}
for some invariant functional $U_{(\boldsymbol{X}_1,\ldots,\boldsymbol{X}%
_n)}(\cdot) \in \widetilde{\mathcal{G}}$, i.e., 
\begin{equation*}
U_{(g(x_1),\ldots,g(x_n))}(\cdot) = U_{(x_1,\ldots,x_n)}(\cdot), \, \forall
g\in\mathcal{G}, \forall (\boldsymbol{x}_1,\ldots,\boldsymbol{x}_n)\in \R%
^{mn}.
\end{equation*}
\end{proposition}

\begin{proposition}
(\cite{Hoff_2012_Preprint}, Theorem 3). A functional \, $U : (\boldsymbol{x}%
_1,\ldots,\boldsymbol{x}_n)\in \R^{mn} \mapsto U_{(\boldsymbol{x}_1,\ldots,%
\boldsymbol{x}_n)}(\cdot)\in \widetilde{\mathcal{G}}$ is invariant iff it is
a function of a maximal invariant functional $\mathcal{U}$.
\end{proposition}

Combining the above results, we obtain

\begin{proposition}
Let $\widetilde{\boldsymbol{\theta}}_0(\boldsymbol{X}_1,\ldots,\boldsymbol{X}%
_n)$ be any equivariant estimate. Then $\widetilde{\boldsymbol{\theta}} (%
\boldsymbol{X}_1,\ldots,\boldsymbol{X}_n)$ is equivariant iff 
\begin{equation*}
\widetilde{\boldsymbol{\theta}}(\boldsymbol{X}_1,\ldots,\boldsymbol{X}_n) = 
\mathcal{H}_{\mathcal{U}_{(\boldsymbol{X}_1,\ldots,\boldsymbol{X}_n)}}(%
\widetilde{\boldsymbol{\theta}}_0),
\end{equation*}
where $\mathcal{H}_{\mathcal{U}_{(\boldsymbol{X}_1,\ldots,\boldsymbol{X}%
_n)}}(\cdot)\in\widetilde{\mathcal{G}}$ is some function of the maximal
invariant functional $\mathcal{U}_{(\boldsymbol{X}_1,\ldots,\boldsymbol{X}%
_n)}.$
\end{proposition}

\begin{remark}
Notice that $\mathcal{H}_{\mathcal{U}_{(\boldsymbol{X}_1,\ldots,\boldsymbol{X%
}_n)}}(\cdot)$ acts additively for additive group, and multiplicatively for
multiplicative group, i.e., 
\begin{equation*}
\mathcal{H}_{\mathcal{U}_{(\boldsymbol{X}_1,\ldots,\boldsymbol{X}_n)}}(%
\widetilde{\boldsymbol{\theta}}_0) = \widetilde{\boldsymbol{\theta}}_0 + 
\mathcal{H}(\mathcal{U}_{(\boldsymbol{X}_1,\ldots, \boldsymbol{X}_n)}),
\end{equation*}
when an additive group is considered, and 
\begin{equation*}
\mathcal{H}_{\mathcal{U}_{(\boldsymbol{X}_1,\ldots,\boldsymbol{X}_n)}}(%
\widetilde{\boldsymbol{\theta}}_0) = \mathcal{H}(\mathcal{U}_{(\boldsymbol{X}%
_1,\ldots,\boldsymbol{X}_n)}) \cdot \widetilde{\boldsymbol{\theta}}_0,
\end{equation*}
for multiplicative group.
\end{remark}

\section{UMRE estimate for additive group}

\noindent Let $\widetilde{\boldsymbol{\theta}}_n := \widetilde{\boldsymbol{%
\theta}}_n(\boldsymbol{X}_1, \ldots, \boldsymbol{X}_n)$ be any one of the
equivariant estimates $\widetilde{\boldsymbol{\theta}}_\varphi$, and assume
that $\mathbb{E}\left({\left\| \widetilde{\boldsymbol{\theta}}_n \right\|}%
^{2} \right)<\infty$. Consider the $L_2$ loss. In view of the above
statements, the UMRE estimate of $\boldsymbol{\theta}_T$ is then given by 
\begin{equation}  \label{MRE Pitman}
\widehat{\boldsymbol{\theta}}_n = \widetilde{\boldsymbol{\theta}}_n - 
\mathbb{E}_{\boldsymbol{0}}\left( \widetilde{\boldsymbol{\theta}}_n \, | \, 
\boldsymbol{T} \right),
\end{equation}
where $\boldsymbol{T} :=(\boldsymbol{X}_2-\boldsymbol{X}_1,\ldots,%
\boldsymbol{X}_n-\boldsymbol{X}_1)$, and $\mathbb{E}_{\boldsymbol{0}}\left( 
\widetilde{\boldsymbol{\theta}}_n \, | \, \boldsymbol{T} \right)$ is the
conditional expectation of $\widetilde{\boldsymbol{\theta}}_n$ given $%
\boldsymbol{T}$, under the assumption that $\boldsymbol{\theta}_T = 
\boldsymbol{0}$. We give in the following an asymptotic approximation to the
conditional expectation 
\begin{equation*}
\mathbb{E}_{\boldsymbol{0}}\left( \widetilde{\boldsymbol{\theta}}_n \, | \, 
\boldsymbol{T}\right),
\end{equation*}
using the result of \cite{Jurekova_Picek_2009_StatDec}. Straightforward
calculs, shows that the score function, of the semiparametric exponential
family (\ref{exponential family}), can be written as 
\begin{eqnarray}
\psi_p(\boldsymbol{x},\boldsymbol{\theta}) & := & \frac{\partial}{\partial 
\boldsymbol{\theta}} \log p(\boldsymbol{x},\boldsymbol{\theta}) = \frac{%
(\partial/\partial\boldsymbol{\theta})p(\boldsymbol{x},\boldsymbol{\theta})}{%
p(\boldsymbol{x},\boldsymbol{\theta})}  \label{score function} \\
& = & \boldsymbol{t}^{\prime }(\boldsymbol{\theta}) \mathbf{f}(\boldsymbol{x}%
,\boldsymbol{\theta}) + \mathbf{f}^{\prime }(\boldsymbol{x},\boldsymbol{%
\theta}) \boldsymbol{t}(\boldsymbol{\theta})  \notag \\
& & - \frac{\int_{\R^m } \left[\boldsymbol{t}^{\prime }(\boldsymbol{\theta}) 
\mathbf{f}(\boldsymbol{x},\boldsymbol{\theta}) + \mathbf{f}^{\prime }(%
\boldsymbol{x},\boldsymbol{\theta}) \boldsymbol{t}(\boldsymbol{\theta})%
\right] \exp\left\{\boldsymbol{t}(\boldsymbol{\theta})^\top \mathbf{f}(%
\boldsymbol{x},\boldsymbol{\theta})\right\} \,dP_{\boldsymbol{X}}(%
\boldsymbol{x}) }{\int_{\R^m} \exp\left\{\boldsymbol{t}(\boldsymbol{\theta}%
)^\top \mathbf{f}(\boldsymbol{x},\boldsymbol{\theta})\right\} \,dP_{%
\boldsymbol{X}}(\boldsymbol{x})},
\end{eqnarray}
where $\mathbf{f}^{\prime }(\boldsymbol{x},\boldsymbol{\theta})$ and $%
\boldsymbol{t}^{\prime }(\boldsymbol{\theta})$ are, respectively, the
derivative w.r.t. $\boldsymbol{\theta}$, of $\mathbf{f}(\boldsymbol{x},%
\boldsymbol{\theta})^\top$ and $\boldsymbol{t}(\boldsymbol{\theta})^\top$
the solution of the system (\ref{exponential family system}). The derivative 
$\boldsymbol{t}^{\prime }(\boldsymbol{\theta})$ can be derived by the
implicit function theorem. Denote $\boldsymbol{h}(\boldsymbol{\theta},%
\boldsymbol{t}) := (h_1(\boldsymbol{\theta},\boldsymbol{t}),\ldots,h_\ell(%
\boldsymbol{\theta},\boldsymbol{t}))^\top$ with 
\begin{equation*}
h_j(\boldsymbol{\theta},\boldsymbol{t}) := \int_{\R^m}f_j(\boldsymbol{x},%
\boldsymbol{\theta}) \exp\left\{\boldsymbol{t}^\top \mathbf{f}(\boldsymbol{x}%
,\boldsymbol{\theta})\right\}\, dP_{\boldsymbol{X}}(\boldsymbol{x}),\,
\forall j=1,\ldots,\ell.
\end{equation*}
Let 
\begin{equation*}
J_{\boldsymbol{h}}^{\boldsymbol{\theta}}(\boldsymbol{\theta},\boldsymbol{t})
:= \frac{\partial}{\partial\boldsymbol{\theta}}\boldsymbol{h}(\boldsymbol{%
\theta},\boldsymbol{t})^\top\quad\text{and} \quad J_{\boldsymbol{h}}^{%
\boldsymbol{t}}(\boldsymbol{\theta},\boldsymbol{t}):= \frac{\partial}{%
\partial\boldsymbol{\theta}}\boldsymbol{h}(\boldsymbol{\theta},\boldsymbol{t}%
)^\top.
\end{equation*}
Then, by the implicit function theorem, we have 
\begin{equation}
\boldsymbol{t}^{\prime }(\boldsymbol{\theta}) := \frac{\partial}{\partial 
\boldsymbol{\theta}} \boldsymbol{t}(\boldsymbol{\theta})^\top = - J_{%
\boldsymbol{h}}^{\boldsymbol{\theta}}(\boldsymbol{\theta},\boldsymbol{t}(%
\boldsymbol{\theta})) \left[ J_{\boldsymbol{h}}^{\boldsymbol{t}}(\boldsymbol{%
\theta},\boldsymbol{t}(\boldsymbol{\theta}))\right]^{-1}.
\end{equation}

Notice that, for true value $\boldsymbol{\theta}=\boldsymbol{\theta}_T$,
since $\boldsymbol{t}(\boldsymbol{\theta}_T)=\boldsymbol{0}$, we obtain for
the true value $\boldsymbol{\theta}_T$ the more simpler expression 
\begin{equation}  \label{fonction score pour theta_T}
\psi_p(\boldsymbol{x},\boldsymbol{\theta}_T) = \boldsymbol{t}^{\prime }(%
\boldsymbol{\theta}_T)\mathbf{f}(\boldsymbol{x},\boldsymbol{\theta}_T) =:
t^{\prime }_1(\boldsymbol{\theta}_T)f_1(\boldsymbol{x},\boldsymbol{\theta}%
_T)+\cdots+t_\ell^{\prime }(\boldsymbol{\theta}_T)f_\ell(\boldsymbol{x}, 
\boldsymbol{\theta}_T) .
\end{equation}
Let 
\begin{equation}  \label{jacobien de psi_p}
J_{\psi}(\boldsymbol{x},\boldsymbol{\theta}) := \left[(\partial^2/\partial
\theta_i\partial\theta_j) \log p(\boldsymbol{x},\boldsymbol{\theta})\right]%
_{i,j=1,\ldots,d}.
\end{equation}
Under some integrability assumptions, by dominated convergence theorem, we
obtain 
\begin{equation}  \label{jacobien de psi_p}
\E \left(J_{\psi}(\boldsymbol{X},\boldsymbol{\theta}_T)\right) = - \E%
\left(\psi_p(\boldsymbol{X},\boldsymbol{\theta}_T)\psi_p(\boldsymbol{X},%
\boldsymbol{\theta}_T)^\top\right) =: - I(\boldsymbol{\theta}_T),
\end{equation}
which is the opposite of the Fisher information matrix.

\begin{theorem}
Under some regularity conditions, we have 
\begin{eqnarray}
\E_0\left(\widetilde{\boldsymbol{\theta}}_\varphi \,|\, \boldsymbol{T}%
\right) & = & - {I(\boldsymbol{\theta}_T)}^{-1} \frac{1}{n}\sum_{i=1}^n
\psi_p(\boldsymbol{X}_i,\widetilde{\boldsymbol{\theta}}_\varphi)
+o_P(n^{-1/2}) \\
& \approx & - {\widehat{I}(\widetilde{\boldsymbol{\theta}}_\varphi)}^{-1} 
\frac{1}{n}\sum_{i=1}^n \widehat{\psi}_p(\boldsymbol{X}_i,\widetilde{%
\boldsymbol{\theta}}_\varphi)
\end{eqnarray}
which gives the following approximation of the UMRE estimate 
\begin{equation}  \label{UMRE approximation}
\widehat{\boldsymbol{\theta}} \approx \widetilde{\boldsymbol{\theta}}%
_\varphi + {\widehat{I}(\widetilde{\boldsymbol{\theta}}_\varphi)}^{-1} \frac{%
1}{n}\sum_{i=1}^n \widehat{\psi}_p(\boldsymbol{X}_i,\widetilde{\boldsymbol{%
\theta}}_\varphi),
\end{equation}
where ${\widehat{I}(\widetilde{\boldsymbol{\theta}}_\varphi)}$ is the
empirical estimate of the Fisher information matrix $I(\boldsymbol{\theta}%
_T) $, given by 
\begin{equation*}
{\widehat{I}(\widetilde{\boldsymbol{\theta}}_\varphi)} := \frac{1}{n}%
\sum_{i=1}^n\widehat{\psi}_p(\boldsymbol{X}_i,\widehat{\boldsymbol{\theta}}%
_\varphi)\widehat{\psi}_p(\boldsymbol{X}_i,\widehat{\boldsymbol{\theta}}_
\varphi)^\top,
\end{equation*}
with $\forall \boldsymbol{x}\in\R^m, \forall \boldsymbol{\theta}\in\Theta,$ 
\begin{eqnarray}
\widehat{\psi}_p(\boldsymbol{x},\boldsymbol{\theta}) & := & \widehat{%
\boldsymbol{t}}^{\prime }(\boldsymbol{\theta}) \mathbf{f}(\boldsymbol{x},%
\boldsymbol{\theta}) + \mathbf{f}^{\prime }(\boldsymbol{x},\boldsymbol{\theta%
}) \widehat{\boldsymbol{t}}(\boldsymbol{\theta})  \notag \\
& & - \frac{\int_{\R^m } \left[\widehat{\boldsymbol{t}}^{\prime }(%
\boldsymbol{\theta}) \mathbf{f}(\boldsymbol{x},\boldsymbol{\theta}) + 
\mathbf{f}^{\prime }(\boldsymbol{x},\boldsymbol{\theta}) \widehat{%
\boldsymbol{t}}(\boldsymbol{\theta})\right] \exp\left\{\widehat{\boldsymbol{t%
}}(\boldsymbol{\theta})^\top \mathbf{f}(\boldsymbol{x},\boldsymbol{\theta}%
)\right\} \,dP_n(\boldsymbol{x}) }{\int_{\R^m} \exp\left\{\widehat{%
\boldsymbol{t}}(\boldsymbol{\theta})^\top \mathbf{f}(\boldsymbol{x},%
\boldsymbol{\theta})\right\} \,dP_n(\boldsymbol{x})},
\end{eqnarray}
$\widehat{\boldsymbol{t}}(\boldsymbol{\theta})$ is the solution of the
empirical version of the system (\ref{exponential family system}), i.e., the
solution in $\boldsymbol{t}$ of 
\begin{equation}  \label{empirical exponential family system}
\int_{\R^m}f_j(\boldsymbol{x},\boldsymbol{\theta}) \exp\left\{\boldsymbol{t}%
^\top \mathbf{f}(\boldsymbol{x},\boldsymbol{\theta})\right\}\, dP_n(%
\boldsymbol{x})=0,\, j=1,\ldots,\ell,
\end{equation}
and $\widehat{\boldsymbol{t}}^{\prime }(\boldsymbol{\theta})$ is the
gradient of $\widehat{\boldsymbol{t}}(\boldsymbol{\theta})$ at the point $%
\boldsymbol{\theta}$ given by 
\begin{equation}
\widehat{\boldsymbol{t}}^{\prime }(\boldsymbol{\theta}) := \frac{\partial}{%
\partial \boldsymbol{\theta}} \widehat{\boldsymbol{t}}(\boldsymbol{\theta}%
)^\top= - \widehat{J_{\boldsymbol{h}}^{\boldsymbol{\theta}}}(\boldsymbol{%
\theta},\widehat{\boldsymbol{t}}(\boldsymbol{\theta})) \left[ \widehat{J_{%
\boldsymbol{h}}^{\boldsymbol{t}}}(\boldsymbol{\theta},\widehat{\boldsymbol{t}%
}(\boldsymbol{\theta}))\right]^{-1},
\end{equation}
where 
\begin{equation*}
\widehat{J_{\boldsymbol{h}}^{\boldsymbol{\theta}}}(\boldsymbol{\theta},%
\boldsymbol{t}) := \frac{\partial}{\partial\boldsymbol{\theta}} \widehat{%
\boldsymbol{h}}(\boldsymbol{\theta},\boldsymbol{t})^\top,\quad \widehat{J_{%
\boldsymbol{h}}^{\boldsymbol{t}}}(\boldsymbol{\theta}, \boldsymbol{t}):= 
\frac{\partial}{\partial\boldsymbol{\theta}}\widehat{\boldsymbol{h}}(%
\boldsymbol{\theta},\boldsymbol{t})^\top, \quad \widehat{\boldsymbol{h}}(%
\boldsymbol{\theta},\boldsymbol{t})^\top:=(\widehat{h_1}(\boldsymbol{\theta},%
\boldsymbol{t}),\ldots, \widehat{h_\ell}(\boldsymbol{\theta},\boldsymbol{t}%
))^\top
\end{equation*}
and 
\begin{equation*}
\widehat{h_j}(\boldsymbol{\theta},\boldsymbol{t}) := \int_{\R^m}f_j(%
\boldsymbol{x},\boldsymbol{\theta}) \exp\left\{\boldsymbol{t}^\top \mathbf{f}%
(\boldsymbol{x},\boldsymbol{\theta})\right\}\, dP_n(\boldsymbol{x}),\,
\forall j=1,\ldots,\ell.
\end{equation*}
\end{theorem}

\section{UMRE estimate for multiplicative group}

\noindent Let $\widetilde{\boldsymbol{\theta }}_{n}:=\widetilde{\boldsymbol{%
\theta }}_{n}(\boldsymbol{X}_{1},\ldots ,\boldsymbol{X}_{n})$ be any one of
the equivariant estimates $\widetilde{\boldsymbol{\theta }}_{\varphi }$ of $%
\boldsymbol{\theta }_{T}$, and assume that $\mathbb{E}\left( {\left\Vert 
\widetilde{\boldsymbol{\theta }}_{n}\right\Vert }^{2}\right) <\infty $.
Consider the $L_{r}$ loss. In view of the above statements, the UMRE
estimate of $\boldsymbol{\theta }_{T}$ is given by 
\begin{equation}
\widehat{\boldsymbol{\theta }}_{n}=\frac{\widetilde{\boldsymbol{\theta }}%
_{n}\cdot \E_{\mathbf{1}}\left( \widetilde{\boldsymbol{\theta }}_{n}\,|\,%
\boldsymbol{T}\right) }{\E_{\mathbf{1}}\left( {\widetilde{\boldsymbol{\theta 
}}_{n}}^{\top }\widetilde{\boldsymbol{\theta }}_{n}\,|\,\boldsymbol{T}%
\right) },  \label{MRE Pitman}
\end{equation}%
where $\boldsymbol{T}=\left( \frac{\boldsymbol{X}_{1}}{|\boldsymbol{X}_{1}|},%
\frac{\boldsymbol{X}_{2}}{\boldsymbol{X}_{1}},\ldots ,\frac{\boldsymbol{X}%
_{n}}{\boldsymbol{X}_{1}}\right) $, and $\mathbb{E}_{\mathbf{1}}\left( \cdot
\,|\,\boldsymbol{T}\right) $ is the conditional expectation given $%
\boldsymbol{T}$, under the assumption that $\boldsymbol{\theta }_{T}=\mathbf{%
1}$.

\section{Simulation results}

\begin{center}
\begin{figure}[tbp]
\label{fig1} \includegraphics[width=10cm]{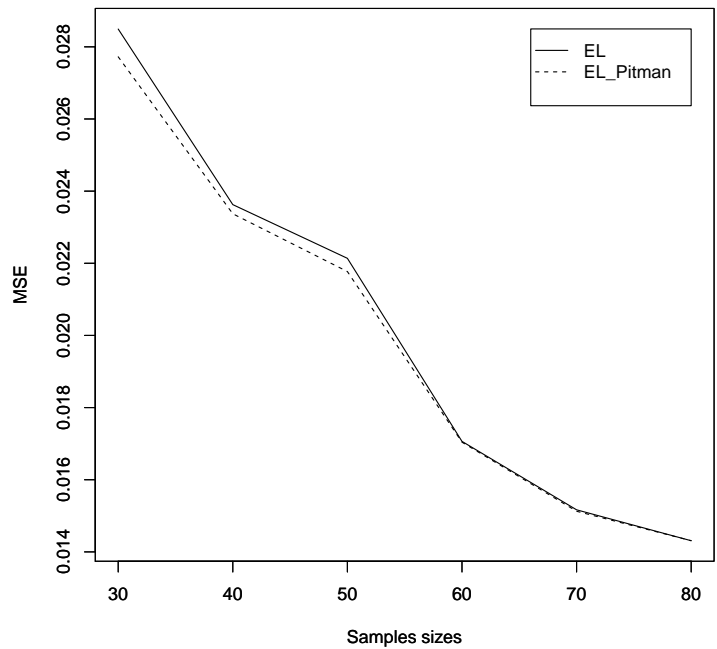}
\caption{MSE v.s. sample sizes}
\end{figure}
\end{center}

\begin{example}
\label{exemple 1} Consider the model 
\begin{equation*}
\mathcal{M}:=\bigcup_{\theta \in \Theta }\mathcal{M}_{\theta
}:=\bigcup_{\theta \in \Theta }\left\{ Q\in M^{1}\text{ such that }\int_{%
\mathbb{R}}\mathbf{f}(x,\theta )\,dQ(x)=\boldsymbol{0}\right\} ,
\end{equation*}%
where $\mathbf{f}:(x,\theta )\in \mathbb{R}\times \Theta \mapsto \mathbf{f}%
(x,\theta ):=(x-\theta ,(x-\theta )^{2}-1)^{\top }$. Let $X$ be a random
variable with distribution $P_{X}:=\mathcal{N}(\theta _{T},1)$ with $\theta
_{T}=1$. The model is invariant to the additive group. We compare the mean
square errors (MSE) of the EL estimate and the proposed UMRE estimate using
the approximation (\ref{UMRE approximation}), for the sample sizes $%
n=30,40,50,60,70,80$, with $1000$ runs. We can see, from figure \ref{fig1},
that the proposed estimate improves the EL one for moderate sample sizes.
\end{example}


\begin{thebibliography}{Jure{\v{c}}kov{\'a} and Picek(2009)}
\bibitem[Broniatowski and Keziou(2006)]{Bronia_Kez2006_STUDIA} Broniatowski,
M. and Keziou, A. (2006). \newblock Minimization of {$\phi$}-divergences on
sets of signed measures. 
\newblock {\em Studia Sci. Math. Hungar.;
arXiv:1003.5457}, \textbf{43}(4), 403--442.

\bibitem[Broniatowski and Keziou(2012)]{Bronia_Keziou2012} Broniatowski, M.
and Keziou, A. (2012). \newblock Divergences and duality for estimation and
test under moment condition models. 
\newblock {\em J. Statist. Plann.
Inference}, \textbf{142}(9), 2554--2573.

\bibitem[Chen and Qin(1993)]{ChenAndQin1993} Chen, J.~H. and Qin, J. (1993). %
\newblock Empirical likelihood estimation for finite populations and the
effective usage of auxiliary information. \newblock {\em Biometrika}, 
\textbf{80}(1), 107--116.

\bibitem[Cressie and Read(1984)]{Cressie-Read1984} Cressie, N. and Read, T.
R.~C. (1984). \newblock Multinomial goodness-of-fit tests. 
\newblock {\em J.
Roy. Statist. Soc. Ser. B}, \textbf{46}(3), 440--464.

\bibitem[Csisz{\'a}r(1963)]{Csiszar1963} Csisz{\'a}r, I. (1963). \newblock %
Eine informationstheoretische {U}ngleichung und ihre {A}nwendung auf den {B}%
eweis der {E}rgodizit\"at von {M}arkoffschen {K}etten. 
\newblock {\em Magyar
Tud. Akad. Mat. Kutat\'o Int. K\"ozl.}, \textbf{8}, 85--108.

\bibitem[Csisz{\'a}r(1967)]{Csiszar1967} Csisz{\'a}r, I. (1967). \newblock %
On topology properties of {$f$}-divergences. 
\newblock {\em Studia Sci.
Math. Hungar.}, \textbf{2}, 329--339.

\bibitem[Godambe and Thompson(1989)]{GodambeAndThompson1989} Godambe, V.~P.
and Thompson, M.~E. (1989). \newblock An extension of quasi-likelihood
estimation. \newblock {\em J. Statist. Plann. Inference}, \textbf{22}(2),
137--172. \newblock With discussion and a reply by the authors.

\bibitem[Haberman(1984)]{Haberman1984} Haberman, S.~J. (1984). \newblock %
Adjustment by minimum discriminant information. \newblock {\em Ann. Statist.}%
, \textbf{12}(3), 971--988.

\bibitem[Hansen \emph{et~al.}(1996)]{Hansen_Healton_Yaron1996} Hansen, L.,
Heaton, J., and Yaron, A. (1996). \newblock Finite-sample properties of some
alternative gmm estimators. 
\newblock {\em Journal of Business and Economic
Statistics}, \textbf{14}, 462--2800.

\bibitem[Hansen(1982)]{Hansen1982} Hansen, L.~P. (1982). \newblock Large
sample properties of generalized method of moments estimators. \newblock
\emph{Econometrica}, \textbf{50}(4), 1029--1054.

\bibitem[Hoff(2012)]{Hoff_2012_Preprint} Hoff, P.~D. (2012). \newblock %
Equivariant estimation. \newblock {\em Preprint}.

\bibitem[Imbens(1997)]{Imbens1997} Imbens, G.~W. (1997). \newblock One-step
estimators for over-identified generalized method of moments models. %
\newblock {\em Rev. Econom. Stud.}, \textbf{64}(3), 359--383.

\bibitem[Jure{\v{c}}kov{\'a} and Picek(2009)]{Jurekova_Picek_2009_StatDec} %
Jure{\v{c}}kov{\'a}, J. and Picek, J. (2009). \newblock Minimum risk
equivariant estimator in linear regression model. 
\newblock {\em Statist.
Decisions}, \textbf{27}(1), 37--54.

\bibitem[Kuk and Mak(1989)]{KukAndMak1989} Kuk, A. Y.~C. and Mak, T.~K.
(1989). \newblock Median estimation in the presence of auxiliary
information. \newblock {\em J. Roy. Statist. Soc. Ser. B}, \textbf{51}(2),
261--269.

\bibitem[Lehmann and Casella(1998)]{Lehmann_Casella1998} Lehmann, E.~L. and
Casella, G. (1998). \newblock {\em Theory of point estimation}. \newblock %
Springer Texts in Statistics. Springer-Verlag, New York, second edition.

\bibitem[Liese and Vajda(1987)]{Liese-Vajda1987} Liese, F. and Vajda, I.
(1987). \newblock {\em Convex statistical distances}, volume~95. \newblock %
BSB B. G. Teubner Verlagsgesellschaft, Leipzig.

\bibitem[McCullagh and Nelder(1983)]{McCullagh_Nelder1983} McCullagh, P. and
Nelder, J.~A. (1983). \newblock {\em Generalized linear models}. \newblock %
Monographs on Statistics and Applied Probability. Chapman \& Hall, London.

\bibitem[Newey and Smith(2004)]{NeweySmith2004} Newey, W.~K. and Smith,
R.~J. (2004). \newblock Higher order properties of {GMM} and generalized
empirical likelihood estimators. \newblock {\em Econometrica}, \textbf{72}%
(1), 219--255.

\bibitem[Owen(1990)]{Owen1990} Owen, A. (1990). \newblock Empirical
likelihood ratio confidence regions. \newblock {\em Ann. Statist.}, \textbf{%
18}(1), 90--120.

\bibitem[Owen(1988)]{Owen1988} Owen, A.~B. (1988). \newblock Empirical
likelihood ratio confidence intervals for a single functional. \newblock
\emph{Biometrika}, \textbf{75}(2), 237--249.

\bibitem[Owen(2001)]{Owen2001} Owen, A.~B. (2001). 
\newblock {\em Empirical
Likelihood}. \newblock Chapman and Hall, New York.

\bibitem[Pardo(2006)]{PardoLeandro2006} Pardo, L. (2006). 
\newblock {\em
Statistical inference based on divergence measures}, volume 185 of \emph{%
Statistics: Textbooks and Monographs}. \newblock Chapman \& Hall/CRC, Boca
Raton, FL.

\bibitem[Qin and Lawless(1994)]{Qin-Lawless1994} Qin, J. and Lawless, J.
(1994). \newblock Empirical likelihood and general estimating equations. %
\newblock {\em Ann. Statist.}, \textbf{22}(1), 300--325.

\bibitem[Rockafellar(1970)]{Rockafellar1970} Rockafellar, R.~T. (1970). %
\newblock {\em Convex analysis}. \newblock Princeton University Press,
Princeton, N.J.

\bibitem[Sheehy(1987)]{Sheehy1987} Sheehy, A. (1987). \newblock \text{%
Kullback-Leibler} constrained estimation of probability measures. \newblock
\emph{Report, Dept. Statistics, Stanford Univ.}

\bibitem[Smith(1997)]{Smith1997} Smith, R.~J. (1997). \newblock Alternative
semi-parametric likelihood approches to generalized method of moments
estimation. \newblock {\em Economic Journal}, \textbf{107}, 503--519.
\end{thebibliography}

\end{document}